\documentclass{article}
\usepackage{amsmath}
\usepackage{amsfonts}
\usepackage{amssymb}
\usepackage{graphicx}
\newcommand{\bpartial}{\mathop{\partial\kern -4pt\raisebox{.8pt}{$|$}}}
\newcommand{\bra}{\mathopen{[\kern-1.6pt[}}
\newcommand{\ket}{\mathclose{]\kern-1.5pt]}}
\newcommand{\bbra}{\mathopen{[\kern-2.2pt[\kern-2.3pt[}}
\newcommand{\bket}{\mathclose{]\kern-2.1pt]\kern-2.3pt]}}

\makeindex
\begin{document}

\title {\large{ \bf Some compatible Poisson structures and integrable bi-Hamiltonian systems on four dimensional and
nilpotent six dimensional symplectic real  Lie groups   }}

\vspace{3mm}
\author { \small{ \bf J. Abedi-Fardad$^1$ }\hspace{-2mm}{ \footnote{ e-mail: j.abedifardad@bonabu.ac.ir}}, \small{ \bf A.
Rezaei-Aghdam$^2$ }\hspace{-1mm}{\footnote{ e-mail: rezaei-a@azaruniv.edu - Corresponding author.}}, \small{ \bf  Gh. Haghighatdoost$^{1,3}$}\hspace{-1mm}{ \footnote{ e-mail: gorbanali@azaruniv.edu}}  \\
{\small{$^{1}$\em Department of Mathematics,
University of Bonab, Tabriz, Iran}}\\
{\small{$^{2}$\em Department of Physics, Azarbaijan Shahid Madani
University, 53714-161, Tabriz, Iran}}\\
 {\small{$^{3}$\em Department
of Mathematics, Azarbaijan Shahid Madani University, 53714-161, Tabriz, Iran}} }
 \maketitle

\begin{abstract}
 We provide an alternative method for obtaining of compatible Poisson structures on Lie groups by means of 
 the adjoint representations
  of Lie algebras. In this way we calculate  some compatible Poisson structures on four dimensional and
nilpotent six dimensional symplectic real  Lie groups. Then using   Magri-Morosi's
 theorem we obtain new bi-Hamiltonian systems with  four dimensional and
nilpotent six dimensional symplectic real  Lie groups as phase spaces.

\end{abstract}

{\bf keywords:}{ Integrable  bi-Hamiltonian system, Compatible Poisson structures, Symplectic Lie group.}
\section{Introduction}
Compatible different Hamiltonian structures
have been proved to be a valuable tool in the construction of infinite hierarchies of
symmetries and conservation laws for PDEs (see e.g. KdV equation $\cite{VCDS}$). Kdv equation can be studied as a
 bi-Hamiltonian system.
 The study of  bi-Hamiltonian systems started with the pioneering
work of Magri $\cite{FM1}$  and developed later in many papers (see for example $\cite{FM2}$, $\cite{FM3}$ and
 $\cite{AB1} $).
The bi-Hamiltonian structure has been observed in many of classical systems and some new interesting  examples of
bi-Hamiltonian systems have
been discovered (see for example $\cite{AVT}$,  $\cite{DAA}  $ and $\cite{AVT1}$).
  In this work, we give a method to  construct
 compatible Poisson structures on a Lie group by means of the adjoint representation of its Lie algebra and construct integrable
 bi-Hamiltonian systems by using Magri-Morosi's theorem $\cite{FMC}$(for a review see $\cite{FM}$). Of course using
 of adjoint representation in the context of the coadjoint orbit method applied previously in the Hamiltonian system
 (see for example $\cite{FM}$ and $\cite{AP}$).  We will give
  a method to produce integral of motion of a non-degenerate bi-Hamiltonian systems, for which the Lie group is the
phase space.

 The outline of the paper is as follows.
In section two after reviewing the definition of compatible Poisson structures we give a method for calculating compatible Poisson structures on a Lie group (in general). Then in section three we have obtained these structures on symplectic four
dimensional real Lie groups. In section four using Magri-Morosi's theorem we have obtained new bi-Hamiltonian systems
on symplectic four dimensional real Lie groups (as phase space). In the same way we have obtained the compatible
 Poisson structures on symplectic nilpotent six dimensional real Lie groups $\cite{VVM}$ (see $\cite{VVM11}$ for a
 rigorous commutative relations) in section five, of course in this section we have obtained vielbeins on symplectic nilpotent
six dimensional real Lie algebras, where the results are summarized in appendix B. Finally in section six some new
 bi-Hamiltonian systems on symplectic nilpotent six
 dimensional real Lie groups have been obtained. Some concluding remarks are given in section seven. The list of four dimensional and nilpotent six dimensional symplectic real
 Lie algebras are given in appendix A and B


\section {\large {\bf   Compatible Poisson structures and bi-Hamiltonian systems on  Lie groups}}
\subsection {\large {\bf   Definitions and Notations}}
For the sake of completeness of the paper let us have a short review on the compatible Poisson structures and
 bi-Hamiltonian system (for a review see \cite{FM}).

{\bf Definition}: \cite{FMC} A  pair of  Poisson brackets $ \lbrace .,. \rbrace$ and   $\lbrace .,. \rbrace ^\prime $ or a pair of Poisson bivectors ${\bf P}$ and ${\bf P^\prime}$ on an $m$ dimensional manifold ${\bf M}$ is called compatible if we have:
\begin{equation}\label{AA1}
[{\bf P}, {\bf P}]=[{\bf P^\prime} ,{\bf P^\prime}]=[{\bf P},{\bf P^\prime}]=0,
\end{equation}
where [.,.] is the Schouten bracket that have the following forms:
\begin{equation}\label{AAA2}
[{\bf P}, {\bf P}]^{\lambda \mu \nu }= {\bf P}^{\rho \lambda} \partial _\rho {\bf P}^{ \mu \nu}
+   {\bf P}^{\rho \nu} \partial _\rho {\bf P}^{ \lambda \mu} +   {\bf P}^{\rho \mu} \partial _\rho {\bf P}^{ \nu \lambda},
\end{equation}
\begin{equation}\label{AAAA2}
[{\bf P^\prime}, {\bf P^\prime}]^{\lambda \mu \nu }= {\bf P^\prime}^{\rho \lambda} \partial _\rho {\bf P^\prime}^{ \mu \nu}
+   {\bf P^\prime}^{\rho \nu} \partial _\rho {\bf P^\prime}^{ \lambda \mu} +   {\bf P^\prime}^{\rho \mu} \partial _\rho {\bf P^\prime}^{ \nu \lambda},
\end{equation}
\begin{equation}\label{AA2}
[{\bf P},{\bf P^\prime}]^{\lambda \mu \nu }= {\bf P}^{\rho \lambda} \partial _\rho {\bf P}^{\prime \mu \nu}
 + {\bf P}^{\prime \rho \lambda} \partial _\rho {\bf P}^{ \mu \nu}  +
  {\bf P}^{\rho \nu} \partial _\rho {\bf P}^{\prime \lambda \mu} + {\bf P}^{\prime \rho \nu} \partial _\rho {\bf P}^{ \lambda \mu}
+  {\bf P}^{\rho \mu} \partial _\rho {\bf P}^{\prime \nu \lambda} + {\bf P}^{\prime \rho \mu} \partial _\rho {\bf P}^{ \nu \lambda},
\end{equation}
where $\partial _\rho = \frac{\partial}{\partial x_\rho}$ such that
$(x_1,..., x_{m})$ is the coordinate of the manifold ${\bf M}$. The Poisson
bracket corresponding to the Poisson bivector ${\bf P}$ has the form
\begin{equation}\label{AA3}
\lbrace f, g\rbrace = {\bf P}^{\mu \nu}\; \partial _{\mu} f \; \partial _{\nu} g.
\end{equation}
 The bracket  $(\ref{AA3})$ satisfies  the Jacobi identity  i.e. $\forall \;f,g,h \in C^\infty ({\bf M})$,
 \begin{equation}\label{AA4}
 \lbrace f, \lbrace g,h\rbrace \rbrace + \lbrace g, \lbrace h,f\rbrace \rbrace +\lbrace h,
  \lbrace f,g \rbrace \rbrace = [{\bf P},{\bf P}]^{\lambda \mu \nu} \partial _\lambda f \; \partial _\mu g\; \partial _\nu h =0
 \end{equation}
  if $[{\bf P},{\bf P}]=0$ and vice versa. A manifold ${\bf M}$ equipped with
 such compatible Poisson structures is called bi-Hamiltonian manifold.
If a dynamical system on the manifold ${\bf M}$ for which the number of functionally independent integrals of motion
  $H_1, ..., H_{n}$ are in bi-involution with respect to this compatible Poisson brackets,
 \begin{equation}
 \lbrace H_i,H_j \rbrace =\lbrace H_i,H_j \rbrace ^\prime =0
 \end{equation}
 then the system is called  bi-Hamiltonian  \cite{FM}. So to introduce the bi-Hamiltonian  structure on the manifold ${\bf M}$ , we must determine a pair of compatible
  and independent Poisson bivectors ${\bf P}$ and ${\bf P^\prime}$.

\subsection {\large {\bf     Compatible Poisson structures on  Lie groups}}
Now we will try to simplify the  relations $(\ref{AAA2})$, $(\ref{AAAA2})$ and $(\ref{AA2})$ when ${\bf M}$ is a Lie group
by using non-coordinate bases.

{\bf Definition}: \cite{MNA} In the coordinate basis, $T_p{\bf M} $ spanned by $\lbrace e_\mu \rbrace =\lbrace \partial_\mu \rbrace $
and  $T^\ast _p{\bf M} $ by $\lbrace dx^\mu \rbrace$, let us consider their linear combinations,
\begin{equation}\label{AA5}
{\bf e} _i= e_i ^{\;\;\mu}\; \partial _\mu\; ,\;\;\;{\bf{\Theta}}^{ i}=e^i_{\;\;\mu} \; dx^\mu \; ,
\;\;\;\;\;\lbrace e_i ^{\;\;\mu} \rbrace \in GL(m,R),
\end{equation}
where $det \; e_i ^{\;\;\mu} >0$. In other words , $\lbrace
{\bf e} _i \rbrace$ is the frame of basis vectors which is
obtained by a $ GL(m,R)$-rotation of the basis  $\lbrace e_\mu \rbrace$ and preserving the orientation. In  the above
 $e^i_{\;\;\mu}$ is inverse of $ e_i ^{\;\;\mu}$ and we have
 \begin{equation}\label{AA6}
 e^i_{\;\;\mu} \;  e_i ^{\;\;\nu} =\delta _\mu ^{\;\;\nu}     \;\;,\;\;\; e^i_{\;\;\mu}\;e_j ^{\;\;\mu} =\delta ^i _{\;\;j}.
 \end{equation}
 The bases $\lbrace {\bf e} _i \rbrace$ and $\lbrace {\bf \Theta} ^{ i}\rbrace$ are called the non-coordinate bases.
 The coefficients $e_i ^{\;\;\mu}$ are called the vielbeins and we have
 \begin{equation}\label{AA7}
 [{\bf e} _i , {\bf e} _j ]= C_{ij}^{\;\;\;k} {\bf e} _k,
 \end{equation}
 where $ C_{ij}^{\;\;\;k}$ is a function of coordinates of the manifold ${\bf M}$. When ${\bf M}$ is  a Lie group
 ${\bf G}$ , these coefficients are the  structure constants of the Lie algebra {\bf g} of the  Lie group
 ${\bf G}$ and we have
 \begin{equation}\label{AA8}
 C_{ij}^{\;\;\;k}=  e^k_{\;\;\nu}\;(e_i ^{\;\;\mu}\; \partial _\mu\; e_j ^{\;\;\nu} - e_j ^{\;\;\mu}\; \partial _\mu \; e_i ^{\;\;\nu}).
 \end{equation}
 Now  we write the Poisson structure ${\bf P}$ in terms of the non-coordinate basis as{\footnote{Here the
 indices $ \mu , \nu,...$ related to the coordinates of the group  ${\bf G}$  and the indices $i,j,...$related to the group
 parameter $\{x_i\}$. Furthermore,  in the following we will consider ${\bf G}$ as a phase space of the dynamical system, hence $m=2n$ must be even.} }
 \begin{equation}\label{AA9}
 {\bf P}^{\mu \nu}= e_i ^{\;\;\mu}\; e_j ^{\;\;\nu} \;P^{ij}.
 \end{equation}
Note that in general $P^{ij}$'s are antisymmetric tensors and functions of the group parameters $x_i$'s. As a first case we consider  $P^{ij}$ and $P^{\prime ij}$  as constant antisymmetric matrices
  \begin{equation}\label{3A111}
P=\left(
\begin{matrix}
0&p_{12}&p_{13}&...&p_{1m}  \cr
 -p_{12}&0&p_{23}&...&p_{2m} \cr
 - p_{13}&-p_{23}&0&...&p_{3m}\cr
 .&.&&&.\cr
 .&.&&&.\cr
 .&.&&&.\cr
 -p_{1m}&-p_{2m}&-p_{3m}&...&0
 \end{matrix}
\right) ,\qquad P^\prime =\left(
\begin{matrix}
0&p^\prime _{12}&p^\prime _{13}&...&p^\prime _{1m}  \cr
 -p^\prime _{12}&0&p^\prime _{23}&...&p^\prime _{2m} \cr
 - p^\prime _{13}&-p^\prime _{23}&0&...&p^\prime _{3m}\cr
 .&.&&&.\cr
 .&.&&&.\cr
 .&.&&&.\cr
 -p^\prime _{1m}&-p^\prime _{2m}&-p^\prime _{3m}&...&0
 \end{matrix}
\right),
\end{equation}
where $p_{ij}$ and $p^\prime _{ij}$ are real constants.
  Now using  $(\ref{AA1})$, $(\ref{AA8})$ and $(\ref{AA9})$ one can rewrite the relations $(\ref{AA1})$ with $(\ref{AAA2})$, $(\ref{AAAA2})$ and $(\ref{AA2})$  as follows:\\
    \begin{equation}\label{AA1011}
 C_{k i}^{\;\;\;s} P^{k z}P^{ i\gamma} +
 C_{k j}^{\;\;\;\gamma} P^{k z}P^{ sj}+
 C_{k i}^{\;\;\;z} P^{k \gamma}P^{ is} =0,
\end{equation}
  \begin{equation}\label{AA1012}
 C_{k i}^{\;\;\;s} P^{\prime k z}P^{\prime i\gamma} +
 C_{k j}^{\;\;\;\gamma} P^{\prime k z}P^{\prime sj}+
 C_{k i}^{\;\;\;z} P^{\prime k \gamma}P^{\prime is} =0,
\end{equation}
  \begin{equation}\label{AA1013}
 C_{k i}^{\;\;\;s}(P^{k z}P^{\prime i\gamma}+P^{\prime k z}P^{ i\gamma}) +
 C_{k j}^{\;\;\;\gamma}(P^{k z}P^{\prime sj}+P^{\prime k z}P^{ sj}) +
 C_{k i}^{\;\;\;z}(P^{k \gamma}P^{\prime is}+P^{\prime k \gamma}P^{ is}) =0.
\end{equation}
then using the adjoint representation of the Lie algebra ${\bf g}$, i.e.:  $({\cal X}_i)_j^{\;k}=-{C_{ij}}^k $ and $({\cal Y}^k)_{ij}=-{C_{ij}}^k$, one can rewrite the above relations in the following matrix forms:{\footnote{Here the upper index $t$
represents the transpose of a matrix.} }

\begin{equation}\label{AA10}
P \;{\cal X}_i \; P^{i\gamma} + P\;{\cal Y}^\gamma \;P + P^{i\gamma}\;{\cal X}_i^{\;t}\;P=0,
\end{equation}
\begin{equation}\label{AAA10}
P^\prime \;{\cal X}_i \; P^{\prime i\gamma} + P^\prime\;{\cal Y}^{\gamma}\;P^\prime + P^{\prime i\gamma}\;{\cal X}_i^{\;t}\;P^\prime =0,
\end{equation}
\begin{equation}\label{AA11}
P \;{\cal X}_i \; P^{\prime i\gamma} + P^\prime \;{\cal X}_i \; P^{i\gamma} + P\;{\cal Y}^\gamma \;P^\prime +
 P^\prime \;{\cal Y}^\gamma \;P + P^{i\gamma }\;{\cal X}_i^{\;t}\;P^\prime + P^{\prime i\gamma}\;{\cal X}_i^{\;t}\;P =0,
\end{equation}
In this way having the structure constants of the Lie algebra ${\bf g}$, one can
solve the matrix equations (\ref{AA10}) - (\ref{AA11}) in order to
obtain $P$ and $P^\prime $. Here we will consider four dimensional real Lie algebras \cite{JP}.    For all
of the  symplectic{\footnote{Note that here we will consider symplectic four dimensional
real Lie algebras \cite{JAF2} and not all of them \cite{JP}, because  we will construct integrable systems
over these related Lie groups. } } four dimensional
real Lie  algebras \cite{GO}, one can see that all the solutions of  (\ref{AA10}) - (\ref{AA11})
are equivalent.

For this reason as a  second case we  consider $P^{ ij}$ as (\ref{3A111}) but   $P^{\prime ij}$ as a linear functions of group
parameters $x_i$ of the Lie group ${\bf G}$ as follows:

 \begin{equation}\label{AAAA11}
P^\prime =\left(
\begin{matrix}
0&p^\prime _{12}+{\sum ^m_{i=1}} a_{2 i} x_i &p^\prime _{13}+{\sum ^m_{i=1}}a_{3 i} x_i&...&p^\prime _{1m}+{\sum ^m_{i=1}}a_{m i} x_i \cr
-p^\prime _{12}-{\sum ^m_{i=1}} a_{2 i} x_i &0&p^\prime _{23}+{\sum ^m_{i=1}}b_{3 i } x_i &...&p^\prime _{2m}+{\sum ^m_{i=1}}b_{m i} x_i\cr
- p^\prime _{13}- {\sum ^m_{i=1}}a_{3 i} x_i &-p^\prime _{23}-{\sum ^m_{i=1}}b_{3 i } x_i  &0&...&p^\prime _{3m}+{\sum^4_{i=1}}c_{m i} x_i \cr
. &. &. & &.\cr
. &. &. & &.\cr
. &. &. & &.\cr
-p^\prime _{1m}-{\sum ^m_{i=1}}a_{m i} x_i &-p^\prime _{2m}-{\sum ^m_{i=1}}b_{m i} x_i &-p^\prime _{3m}-{\sum^4_{i=1}}c_{m i} x_i  &...&0
 \end{matrix}
\right),
\end{equation}
where $ p^{\prime}_{ij}$ and $a_{ij}$'s are real constants; then   relations $(\ref{AA1})$ with $(\ref{AAAA2})$ and $(\ref{AA2})$  have the following matrix forms:

\begin{equation}\label{AA12}
P^\prime \;{\cal X}_i \; P^{\prime i\gamma} + P^\prime\;{\cal Y}^\gamma \;P^\prime + P^{\prime i\gamma}\;{\cal X}_i^{\;t}\;P^\prime
+ (e^t P^\prime )^{k\gamma}\partial _k P^\prime +A+B=0,
\end{equation}
\begin{equation}\label{AA13}
P \;{\cal X}_i \; P^{\prime i\gamma} + P^\prime \;{\cal X}_i \; P^{i\gamma} + P\;{\cal Y}^\gamma \;P^\prime +
 P^\prime \;{\cal Y}^\gamma\;P + P^{i\gamma}\;{\cal X}_i^{\;t}\;P^\prime + P^{\prime i\gamma}\;{\cal X}_i^{\;t}\;P +
 (e^t P)^{k\gamma }\partial _k P^\prime +A^\prime +B^\prime=0,
\end{equation}
where $e^t$ is a transpose of the vielbein $e_\alpha ^{\;\mu }$ and  $A,B,A^\prime$  and $ B^\prime$ have the following forms:
\begin{equation}
A=\left(
\begin{matrix}
 (e^t P^\prime )^{k1}\partial _k P^{\prime 1\gamma}&(e^t P^\prime )^{k1}\partial _k P^{\prime 2\gamma}&...& (e^t P^\prime )^{k1}\partial _k P^{\prime m\gamma}\cr
 (e^t P^\prime )^{k2}\partial _k P^{\prime 1\gamma}&.&&. \cr
.&.&&. \cr
.&.&&. \cr
 (e^t P^\prime )^{km}\partial _k P^{\prime 1\gamma}&(e^t P^\prime )^{km}\partial _k P^{\prime 2\gamma}&...& (e^t P^\prime )^{km}\partial _k P^{\prime m\gamma}
 \end{matrix}
\right)
\end{equation}

\begin{equation}
B=\left(
\begin{matrix}
 (e^t P^\prime )^{k1}\partial _k P^{\prime \gamma 1}&(e^t P^\prime )^{k2}\partial _k P^{\prime \gamma 1}&...& (e^t P^\prime )^{km}\partial _k P^{\prime \gamma 1}\cr
 (e^t P^\prime )^{k1}\partial _k P^{\prime \gamma 2}&.&&. \cr
.&.&&. \cr
.&.&&. \cr
 (e^t P^\prime )^{k1}\partial _k P^{\prime \gamma m}&(e^t P^\prime )^{k2}\partial _k P^{\prime \gamma m}&...& (e^t P^\prime )^{km}\partial _k P^{\prime \gamma m}
 \end{matrix}
\right)
\end{equation}

\begin{equation}
A^\prime =\left(
\begin{matrix}
 (e^t P )^{k1}\partial _k P^{\prime 1\gamma}&(e^t P)^{k1}\partial _k P^{\prime 2\gamma}&...& (e^t P )^{k1}\partial _k P^{\prime m\gamma}\cr
 (e^t P )^{k2}\partial _k P^{\prime 1\gamma}&.&&. \cr
.&.&&. \cr
.&.&&. \cr
 (e^t P )^{km}\partial _k P^{\prime 1\gamma}&(e^t P )^{km}\partial _k P^{\prime 2\gamma}&...& (e^t P )^{km}\partial _k P^{\prime m\gamma}
 \end{matrix}
\right)
\end{equation}

\begin{equation}
B^\prime =\left(
\begin{matrix}
 (e^t P )^{k1}\partial _k P^{\prime \gamma 1}&(e^t P )^{k2}\partial _k P^{\prime \gamma 1}&...& (e^t P )^{km}\partial _k P^{\prime \gamma 1}\cr
 (e^t P)^{k1}\partial _k P^{\prime \gamma 2}&.&&. \cr
.&.&&. \cr
.&.&&. \cr
 (e^t P )^{k1}\partial _k P^{\prime \gamma m}&(e^t P )^{k2}\partial _k P^{\prime \gamma m}&...& (e^t P )^{km}\partial _k P^{\prime \gamma m}
 \end{matrix}
\right)
\end{equation}
Now we will try to solve (\ref{AA10}), (\ref{AA12})  and (\ref{AA13}) for four dimensional real Lie groups.
\section {\large {\bf Some compatible Poisson structures on symplectic four dimensional real Lie groups}}

Having the structure constants of the Lie algebra ${\bf g}$, we will solve matrix equations (\ref{AA10}), (\ref{AA12})  and (\ref{AA13}) in order to obtain $P$ (\ref{3A111}) and $P^\prime $ (\ref{AAAA11}). For completeness the list of symplectic four dimensional real
Lie algebras \cite{JAF2} are brought in appendix A.

Let us consider an example; for Lie algebra $A_{4,1}$ we have the
following non zero commutators and the matrices ${\cal X}_i$ and ${\cal Y}^i$:
\begin{equation}
[{\bf e}_2,{\bf e}_4]={\bf e}_1\;,\;\;[{\bf e}_3,{\bf e}_4]={\bf e}_2\;,
\end{equation}
\begin{equation}
{\cal X}_1=\left(
\begin{matrix}
0&0&0&0\cr
0&0&0&0\cr
0&0&0&0\cr
0&0&0&0
 \end{matrix}
\right) ,\; {\cal X}_2=\left(
\begin{matrix}
0&0&0&0\cr
0&0&0&0\cr
0&0&0&0\cr
-1&0&0&0
 \end{matrix}
\right) ,\; {\cal X}_3=\left(
\begin{matrix}
0&0&0&0\cr
0&0&0&0\cr
0&0&0&0\cr
0&-1&0&0
 \end{matrix}
\right) ,\; {\cal X}_4=\left(
\begin{matrix}
0&0&0&0\cr
1&0&0&0\cr
0&1&0&0\cr
0&0&0&0
 \end{matrix}
\right),
\end{equation}
\begin{equation}
{\cal Y}^1=\left(
\begin{matrix}
0&0&0&0\cr
0&0&0&-1\cr
0&0&0&0\cr
0&1&0&0
 \end{matrix}
\right) ,\; {\cal Y}^2=\left(
\begin{matrix}
0&0&0&0\cr
0&0&0&0\cr
0&0&0&-1\cr
0&0&1&0
 \end{matrix}
\right) ,\; {\cal Y}^3=\left(
\begin{matrix}
0&0&0&0\cr
0&0&0&0\cr
0&0&0&0\cr
0&0&0&0
 \end{matrix}
\right) ,\; {\cal Y}^4=\left(
\begin{matrix}
0&0&0&0\cr
0&0&0&0\cr
0&0&0&0\cr
0&0&0&0
 \end{matrix}
\right),
\end{equation}
also according to $\cite{bm}$ for  Lie algebra $A_{4,1}$ the matrix $e_\alpha ^{\;\mu} $ has the following form:
\begin{equation}
(e_\alpha ^{\;\mu})=
\left(
\begin{matrix}
1&0&0&0\cr
x_4&1&0&0\cr
x_4^2/2&x_4&1&0\cr
0&0&0&1
 \end{matrix}
\right),
\end{equation}
inserting ${\cal X}_i$, ${\cal Y}^i$ and $e_\alpha ^{\;\mu} $  in (\ref{AA10}), (\ref{AA12})  and (\ref{AA13})  one can obtain
the compatible Poisson structure $P$ and $P^\prime$ for the Lie algebra $A_{4,1}$. One of the solutions has the following forms:
\begin{equation}
P=\left(
\begin{matrix}
0&p _{12} &0&p _{14} \cr
\ast &0&p _{23}&0\cr
 \ast &\ast &0&0\cr
\ast &\ast &\ast &0
 \end{matrix}
\right),\;P^\prime=\left(
\begin{matrix}
0&p^\prime _{12}-(a_{44}+p^\prime _{24}) x_2+a_{23} x_3+a_{24} x_4 &p^\prime _{13}+a_{44} x_3&a_{44} x_4 \cr
\ast &0&0 &p^\prime _{24}\cr
 \ast &\ast &0&0\cr
\ast &\ast &\ast &0
 \end{matrix}
\right).
\end{equation}
Here $p _{12},p _{14},p _{23},p^\prime _{12},p^\prime _{13},p^\prime _{24},a_{23},a_{44}$ and $a_{24}$ arbitrary real
constants.
In this way we have obtained some compatible Poisson structure ($P$ is constant and $P^{\prime}$ is linear
 function of Lie group
coordinates) on symplectic   four dimensional real Lie algebras, the results being summarized in Table 1. Note that all parameters
$a _{ij} , p _{ij}$ and $p^\prime _{ij}$ are arbitrary real constants.
\newpage

{\small {\bf Table 1}}: {\small
 Compatible Poisson structure on symplectic four dimensional real Lie algebras.}\\
    \begin{tabular}{|l|l|l|l|  p{0.15mm} }
    \hline\hline
{\footnotesize ${\bf g}$ }&{\footnotesize $P$ }&{\footnotesize $P^\prime $ }&{\footnotesize Comments}\\ \hline

{\footnotesize$A_{4,1}$}&{\footnotesize $ \left(
\begin{matrix}
0&p _{12} &0&p _{14} \cr
\ast &0&p _{23}&0\cr
 \ast &\ast &0&0\cr
\ast &\ast &\ast &0
 \end{matrix}
\right)$} &{\footnotesize $\left(
\begin{matrix}
0&p^\prime _{12}-(a_{44}+p^\prime _{24}) x_2+a_{23} x_3+a_{24} x_4 &p^\prime _{13}+a_{44} x_3&a_{44} x_4 \cr
\ast &0&0 &p^\prime _{24}\cr
 \ast &\ast &0&0\cr
\ast &\ast &\ast &0
 \end{matrix}
\right)$}&{\footnotesize $p _{14} p _{23}\neq 0$}\\
\hline
{\footnotesize $A_{4,2}^{-1}$}&{\footnotesize $\left(
\begin{matrix}
0&p _{12} &p_{13}&0 \cr
\ast &0&p _{23}&p_{24}\cr
 \ast &\ast &0&0\cr
\ast &\ast &\ast &0
 \end{matrix}
\right)$} &{\footnotesize $\left(
\begin{matrix}
0&p^\prime _{12}+a_{23} x_3+a_{24} x_4 &\frac{a_{53} p_{13}}{p_{23}} x_3&0 \cr
\ast &0&a_{53} x_3 &p^\prime _{24}+a_{64} x_4\cr
 \ast &\ast &0&0\cr
\ast &\ast &\ast &0
 \end{matrix}
\right) $ }&{\footnotesize $p _{13} p _{24}\neq 0$}\\
\hline
{\footnotesize $A_{4,3}$}&{\footnotesize $ \left(
\begin{matrix}
0&p _{12} &0&p_{14} \cr
\ast &0&p _{23}&0\cr
 \ast &\ast &0&0\cr
\ast &\ast &\ast &0
 \end{matrix}
\right) $} &{\footnotesize $ \left(
\begin{matrix}
0&p^\prime _{12}+a_{23} x_3+a_{24} x_4 &a_{33} x_3&a_{33} x_4 \cr
\ast &0&a_{53} x_3 + \frac{a_{33} p_{23}}{p_{14}} x_4 &0\cr
 \ast &\ast &0&0\cr
\ast &\ast &\ast &0
 \end{matrix}
\right) $ }&{\footnotesize $p _{14} p _{23}\neq 0$}\\
\hline
{\footnotesize $A_{4,5}^{a,-a}$}&{\footnotesize $ \left(
\begin{matrix}
0&0 &p _{13}&p_{14} \cr
\ast &0&p _{23}&0\cr
 \ast &\ast &0&0\cr
\ast &\ast &\ast &0
 \end{matrix}
\right)$} &{\footnotesize $\left(
\begin{matrix}
0&0 &p^\prime _{13}+a_{32} x_2+a_{34} x_4 &a_{44} x_4 \cr
\ast &0&p^\prime_{23}+a_{52} x_2  &0\cr
 \ast &\ast &0&0\cr
\ast &\ast &\ast &0
 \end{matrix}
\right) $ }&{\footnotesize $p _{14} p _{23}\neq 0$}\\
\hline
{\footnotesize $A_{4,5}^{-1,-1}$}&{\footnotesize $ \left(
\begin{matrix}
0&p _{12} &p _{13}&0 \cr
\ast &0&p _{23}&p _{24}\cr
 \ast &\ast &0&0\cr
\ast &\ast &\ast &0
 \end{matrix}
\right)$ }&{\footnotesize $\left(
\begin{matrix}
0&p^\prime _{12}+a_{23} x_3+a_{24} x_4&\frac{a_{53} p_{13} x_3}{p_{23}} &0 \cr
\ast &0&a_{53} x_3 &a_{64} x_4\cr
 \ast &\ast &0&0\cr
\ast &\ast &\ast &0
 \end{matrix}
\right) $ }&{\footnotesize $p _{13} p _{24}\neq 0$}\\
\hline
{\footnotesize $A_{4,6}^{a,0}$}&{\footnotesize $ \left(
\begin{matrix}
0&0 &0&p _{14} \cr
\ast &0&p _{23}&0\cr
 \ast &\ast &0&0\cr
\ast &\ast &\ast &0
 \end{matrix}
\right)$ }&{\footnotesize $\left(
\begin{matrix}
0&0 &0 &a_{41} x_1+a_{44} x_4 \cr
\ast &0&a_{52} x_2+a_{53} x_3 &0\cr
 \ast &\ast &0&0\cr
\ast &\ast &\ast &0
 \end{matrix}
\right)$ }&{\footnotesize $p _{14} p _{23}\neq 0$}\\
\hline
{\footnotesize $A_{4,9}^{0}$}&{\footnotesize $ \left(
\begin{matrix}
0&0&p_{13}&p _{23} \cr
\ast &0&p _{23}&0\cr
 \ast &\ast &0&0\cr
\ast &\ast &\ast &0
 \end{matrix}
\right)$} &{\footnotesize $\left(
\begin{matrix}
0&0 &p^\prime _{13}+a_{52} x_1+a_{32}x_2+a_{54}x_3+a_{34}x_4 &p^\prime _{23}+a_{52} x_2+a_{54} x_4 \cr
\ast &0&p^\prime _{23}+a_{52} x_2+a_{54} x_4 &0\cr
 \ast &\ast &0&0\cr
\ast &\ast &\ast &0
 \end{matrix}
\right) $ }&{\footnotesize $ p _{23}\neq 0$}\\

\hline
{\footnotesize $A_{4,12}$}&{\footnotesize $ \left(
\begin{matrix}
0&0 &0 &p _{14} \cr
\ast &0&-p _{14}&0\cr
 \ast &\ast &0&0\cr
\ast &\ast &\ast &0
 \end{matrix}
\right)$} &{\footnotesize $\left(
\begin{matrix}
0&0 &-a_{54} x_3+a_{64}x_4 &-a_{64}x_3-a_{54}x_4 \cr
\ast &0&a_{64}x_3+a_{54}x_4  &-a_{54} x_3+a_{64}x_4\cr
 \ast &\ast &0&0\cr
\ast &\ast &\ast &0
 \end{matrix}
\right) $ }&{\footnotesize $p _{14}\neq 0$}\\
\hline
{\footnotesize $II\oplus R$}&{\footnotesize $\left(
\begin{matrix}
0&p _{12} &p _{13}  &0 \cr
\ast &0&0&0 \cr
 \ast &\ast &0&p _{34} \cr
\ast &\ast &\ast &0
 \end{matrix}
\right)$} &{\footnotesize $\left(
\begin{matrix}
0&a_{44} x_2 &p^\prime_{13}+a_{32} x_2-(a_{44}+p^\prime _{23})x_3 &a_{44}x_4 \cr
\ast &0&p^\prime _{23} &0\cr
 \ast &\ast &0&0\cr
\ast &\ast &\ast &0
 \end{matrix}
\right) $ }&{\footnotesize $p _{12} p _{34}\neq 0$}\\

\hline
{\footnotesize $III\oplus R$}&{\footnotesize $\left(
\begin{matrix}
0&p _{13} &p _{13}  &0 \cr
\ast &0&0&p _{24} \cr
 \ast &\ast &0&0 \cr
\ast &\ast &\ast &0
 \end{matrix}
\right)$} &{\footnotesize $\left(
\begin{matrix}
0&a_{33} x_3 &a_{33} x_3&0\cr
\ast &0&0 &a_{64} x_4\cr
 \ast &\ast &0&0\cr
\ast &\ast &\ast &0
 \end{matrix}
\right) $ }&{\footnotesize $p _{13} p _{24}\neq 0$}\\

\hline
{\footnotesize $VI_0 \oplus R$}&{\footnotesize $\left(
\begin{matrix}
0&p _{12} &0  &0 \cr
\ast &0&0&0 \cr
 \ast &\ast &0&p _{34} \cr
\ast &\ast &\ast &0
 \end{matrix}
\right)$} &{\footnotesize $\left(
\begin{matrix}
0&p^\prime_{12}+a_{22} x_1+a_{22} x_2 &0 &0 \cr
\ast &0&0  &0\cr
 \ast &\ast &0&p^\prime_{34}+a_{73} x_3+a_{74} x_4 \cr
\ast &\ast &\ast &0
 \end{matrix}
\right) $ }&{\footnotesize $p _{12} p _{34}\neq 0$}\\

\hline
{\footnotesize $VII_0 \oplus R$}&{\footnotesize $ \left(
\begin{matrix}
0&p _{12} &0  &0 \cr
\ast &0&0&0 \cr
 \ast &\ast &0&p _{34} \cr
\ast &\ast &\ast &0
 \end{matrix}
\right)$} &{\footnotesize $\left(
\begin{matrix}
0&a_{21} x_1+a_{22} x_2 &0 &0 \cr
\ast &0&0  &0\cr
 \ast &\ast &0&a_{73} x_3+a_{74} x_4 \cr
\ast &\ast &\ast &0
 \end{matrix}
\right) $ }&{\footnotesize $p _{12} p _{34}\neq 0$}\\

\hline
{\footnotesize $A_2 \oplus A_2$}&{\footnotesize $ \left(
\begin{matrix}
0&p _{12} &0  &0 \cr
\ast &0&0&p _{24} \cr
 \ast &\ast &0&p _{34} \cr
\ast &\ast &\ast &0
 \end{matrix}
\right)$} &{\footnotesize $\left(
\begin{matrix}
0&a_{21} x_1&0 &0 \cr
\ast &0&0  & p^\prime_{24}+a_{61} x_1+a_{63} x_3 \cr
 \ast &\ast &0&a_{73} x_3\cr
\ast &\ast &\ast &0
 \end{matrix}
\right) $ }&{\footnotesize $p _{12} p _{34}\neq 0$}\\

\hline
\end{tabular}

\section {\large {\bf Integrable bi-Hamiltonian systems on symplectic four dimensional real Lie groups}}

In this section, we will construct the integrable bi-Hamiltonian systems
with  four dimensional real Lie groups as phase space. For this purpose,
in the previous section we consider  four dimensional real Lie groups such that they have symplectic structures
 \cite{GO}, \cite{JAF2}. Here, we will construct the models
 on these Lie groups separately as follows.\\
An important  class of bi-Hamiltonian manifold occurs when one
 of the compatible Poisson structures is invertible i.e., the Poisson
bracket $\lbrace .,.\rbrace$ associated with ${\bf P}$ is
invertible. Then one can define a linear map ${\bf N}:TM
\longrightarrow TM$  acting on the tangent bundle by \cite{FM}
\begin{equation}\label{A3130}
{\bf N}={\bf P^\prime \; P^{-1}}.
\end{equation}
{\bf Theorem (Magri-Morosi)}: \cite{FMC}, \cite{FM2} \textit{ A remarkable  consequence of the compatibility of  ${\bf P}$ and ${\bf P^\prime}$ is that the  torsion of Nijenhuis tensor
 ${\bf N}$
 \begin{equation}
 T_{{\bf N}}(X,Y)=[{\bf N} X,{\bf N}Y]-{\bf N}[{\bf N}X,Y]-{\bf N}[X,{\bf N}Y]+{\bf N}^2[X,Y]
 \end{equation}
 identically vanishes ; where $X$ and $Y$ are arbitrary vector fields and the bracket $[X, Y ]$ denotes the Lie bracket (commutator).
 One of the main properties of ${\bf N}$ is that the normalized traces of the powers of ${\bf N}$ are integrals of motion
 \begin{equation}\label{A31}
 H_k=\frac{1}{2k} Tr {\bf N}^{k},
 \end{equation}
and satisfying Lenard-Magri recurrent relations \cite{FM}.
 \begin{equation}\label{A313}
{\bf P^\prime}\; dH_i={\bf P}\;dH_{i+1},
 \end{equation}
 }
In  Table 1, the matrices $P$ and $P^\prime$ have been given for the various symplectic Lie algebras.
 Now inserting $P$ and $P^\prime$ in $(\ref{AA9})$
and using the related vielbeins  \cite{bm} the compatible Poisson structures  ${\bf P}$ and ${\bf P^\prime}$
on Lie groups are obtained, then using  $(\ref{A3130})$ and $(\ref{A31})$ one can find the  Hamiltonian and integrals of motions
of bi-Hamiltonian systems. In the following we  will perform this work separately   for symplectic four dimensional real Lie groups.
In this way we  obtain new bi-Hamiltonian  systems over four dimensional real Lie groups as phase spaces.
\vspace*{0.5cm}

{\bf Lie group ${\mathbf A_{4,1}}$:}\\
Inserting $P$ and $P^\prime$ in $(\ref{AA9})$ one can obtain the compatible Poisson structures
${\bf P}$ and ${\bf P^\prime}$ on the Lie group ${\mathbf A_{4,1}}$ as follows:
\begin{equation}
{\bf P}=\left(
\begin{matrix}
 0&p_{12}+\frac{p_{23} x_4^2}{2}&p_{23} x_4& p_{14} \cr
\ast &0&p_{23} &0\cr
\ast &\ast &0&0\cr
\ast &\ast &\ast &0
 \end{matrix}
\right),
 \end{equation}
 \begin{equation}
  {\bf P^\prime} =\left(
\begin{matrix}
0&p^\prime _{12}-( a_{44}  + p^\prime _{24})x_2+a_{23}x_3+a_{24}x_4+p^\prime _{13} x_4+a_{44}x_3 x_4&p^\prime _{13}+a_{44}x_3 &(p^\prime _{24} +a_{44})x_4 \cr
\ast &0 &0& p^\prime _{24}\cr
\ast &\ast &0&0\cr
\ast &\ast &\ast & 0
 \end{matrix}
\right).
\end{equation}
Now by means of $(\ref{A3130})$  and $(\ref{A31})$, the integrals of motion can be found for this Lie group as follows:
\begin{equation}
H_1=\frac{a_{44} x_4}{p_{14}} \;,\;\;H_2=\frac{2 p_{14}p^\prime _{24}(p^\prime _{13}+a_{44} x_3)+a_{44}^2 p_{23} x_4^2}{2p_{14}^2 p_{23}}.
\end{equation}

{\bf Lie group ${\mathbf A_{4,2}^{-1}}$:}\\
Inserting $P$ and $P^\prime$ in $(\ref{AA9})$ one can obtain the compatible Poisson structures
${\bf P}$ and ${\bf P^\prime}$ on the Lie group ${\mathbf A_{4,2}^{-1}}$ as follows:
\begin{equation}
{\bf P}=\left(
\begin{matrix}
 0&p_{12}+p_{13} x_4&p_{13} & 0 \cr
\ast &0&e^{2 x_4}p_{23} &e^{ x_4}p_{24} \cr
\ast &\ast &0&0\cr
\ast &\ast &\ast &0
 \end{matrix}
\right),
 \end{equation}
 \begin{equation}
  {\bf P^\prime} =\left(
\begin{matrix}
0&p^\prime _{12}+a_{23}x_3+a_{24}x_4 +\frac{a_{53} p_{13}x_3 x_4}{p_{23}}&\frac{a_{53} p_{13}x_3 }{p_{23}} &0 \cr
\ast &0 & a_{53} e^{2x_4} x_3&e^{x_4}( p^\prime _{24} +a_{64}x_4) \cr
\ast &\ast &0&0\cr
\ast &\ast &\ast & 0
 \end{matrix}
\right).
\end{equation}
By means of $(\ref{A3130})$  and (\ref{A31}), the integrals of motion can be found for this Lie group as follows:
\begin{equation}
H_1=\frac{p_{23} p^\prime _{24}+a_{53}p_{24}x_3+a_{64} p_{23} x_4}{p_{23} p_{24}} \;,\;\;H_2=\frac{1}{2}
( \frac{ a_{53}^2 x_3^2}{p_{23}^2}+ \frac{(p^\prime _{24}+a_{64} x_4)^2}{p_{24}^2}).
\end{equation}

{\bf Lie group ${\mathbf A_{4,3}}$:}\\
Inserting $P$ and $P^\prime$ in $(\ref{AA9})$ one can obtain the compatible Poisson structures
${\bf P}$ and ${\bf P^\prime}$ on the Lie group ${\mathbf A_{4,3}}$ as follows:
\begin{equation}
{\bf P}=\left(
\begin{matrix}
 0&p_{12} e^{x_4}&0& p_{14} e^{x_4}\cr
\ast &0&p_{23} &0 \cr
\ast &\ast &0&0\cr
\ast &\ast &\ast &0
 \end{matrix}
\right),
 \end{equation}
 \begin{equation}
  {\bf P^\prime} =\left(
\begin{matrix}
0&e^{x_4}(p^\prime _{12}+a_{23}x_3+a_{24}x_4+ a_{33}x_3 x_4)&a_{33} e^{x_4} x_3 &a_{33}e^{x_4} x_4 \cr
\ast &0 & a_{53}  x_3 +\frac{a_{33} p_{23} x_4}{p_{14}}&0\cr
\ast &\ast &0&0\cr
\ast &\ast &\ast & 0
 \end{matrix}
\right).
\end{equation}
Now by means of $(\ref{A3130})$  and (\ref{A31}), the integrals of motion can be found for this Lie group as follows:
\begin{equation}
H_1=\frac{a_{53}x_3}{p_{23}}+ \frac{2 a_{33} x_4}{p_{14}} \;,\;\;H_2=\frac{1}{2}
( \frac{ a_{33}^2 x_4^2}{p_{14}^2}+ (\frac{a_{53}x_3}{p_{23}}+ \frac{ a_{33} x_4}{p_{14}} )^2).
\end{equation}

{\bf Lie group ${\mathbf A_{4,5}^{a,-a}}$:}\\
Inserting $P$ and $P^\prime$ in $(\ref{AA9})$ one can obtain the compatible Poisson structures
${\bf P}$ and ${\bf P^\prime}$ on the Lie group ${\mathbf A_{4,5}^{a,-a}}$ as follows:
\begin{equation}
{\bf P}=\left(
\begin{matrix}
 0&0& p_{13} e^{(1-a)x_4} & p_{14} e^{x_4}\cr
\ast &0&p_{23} &0 \cr
\ast &\ast &0&0\cr
\ast &\ast &\ast &0
 \end{matrix}
\right),
 \end{equation}
 \begin{equation}
  {\bf P^\prime} =\left(
\begin{matrix}
0&0&e^{(1-a)x_4}(p^\prime _{13}+a_{32}x_2+a_{34}x_4) &a_{44}e^{x_4} x_4 \cr
\ast &0 & p^\prime _{23}+a_{52} x_2 &0\cr
\ast &\ast &0&0\cr
\ast &\ast &\ast & 0
 \end{matrix}
\right).
\end{equation}
Now by means of $(\ref{A3130})$  and (\ref{A31}), the integrals of motion can be found for this Lie group as follows:
\begin{equation}
H_1=\frac{p^\prime _{23} + a_{52} x_2}{p_{23}}+ \frac{ a_{44} x_4}{p_{14}} \;,\;\;H_2=\frac{1}{2}
( \frac{(p^\prime _{23}+a_{52} x_2)^2}{p_{23}^2} + \frac{a_{44}^2 x_2^2}{p_{14}^2}).
\end{equation}

{\bf Lie group ${\mathbf A_{4,5}^{-1,-1}}$:}\\
Inserting $P$ and $P^\prime$ in $(\ref{AA9})$ one can obtain the compatible Poisson structures
${\bf P}$ and ${\bf P^\prime}$ on the Lie group ${\mathbf A_{4,5}^{-1,-1}}$ as follows:
\begin{equation}
{\bf P}=\left(
\begin{matrix}
 0&p_{12}& p_{13}  & 0\cr
\ast &0&p_{23} e^{-2x_4} &p_{24} e^{-x_4} \cr
\ast &\ast &0&0\cr
\ast &\ast &\ast &0
 \end{matrix}
\right),
 \end{equation}
 \begin{equation}
  {\bf P^\prime} =\left(
\begin{matrix}
0&p^\prime _{12}+a_{23}x_3+a_{24}x_4 &\frac{a_{53} p_{13}x_3}{p_{23}} &0\cr
\ast &0 & a_{53}e^{-2x_4} x_3 &a _{64}e^{-x_4} x_4  \cr
\ast &\ast &0&0\cr
\ast &\ast &\ast & 0
 \end{matrix}
\right).
\end{equation}
Now by means of $(\ref{A3130})$  and (\ref{A31}), the integrals of motion can be found for this Lie group as follows:
\begin{equation}
H_1=\frac{ a_{53} x_3 }{p_{23}}+ \frac{ a_{64} x_4}{p_{24}} \;,\;\;H_2=\frac{1}{2}
((\frac{ a_{53} x_3 }{p_{23}})^2+( \frac{ a_{64} x_4}{p_{24}})^2 ).
\end{equation}

{\bf Lie group ${\mathbf A_{4,6}^{a,0}}$:}\\
Inserting $P$ and $P^\prime$ in $(\ref{AA9})$ one can obtain the compatible Poisson structures
${\bf P}$ and ${\bf P^\prime}$ on the Lie group ${\mathbf A_{4,6}^{a,0}} $ as follows:
\begin{equation}
{\bf P}=\left(
\begin{matrix}
 0&0& 0 & p_{14} e^{a x_4}\cr
\ast &0& p_{23} &0 \cr
\ast &\ast &0&0\cr
\ast &\ast &\ast &0
 \end{matrix}
\right),
 \end{equation}
 \begin{equation}
  {\bf P^\prime} =\left(
\begin{matrix}
0&0&0 &e^{ax_4}(a_{41}x_1+a_{44}x_4)\cr
\ast &0 & a_{52}x_2+a_{53}x_3 &0 \cr
\ast &\ast &0&0\cr
\ast &\ast &\ast & 0
 \end{matrix}
\right).
\end{equation}
Now by means of $(\ref{A3130})$  and (\ref{A31}), the  integrals of motion can be found for this Lie group as follows:
\begin{equation}
H_1=\frac{a_{41} x_1 +a_{44}x_4}{p_{14}}+ \frac{ a_{52} x_2+a_{53} x_3}{p_{23}} \;,\;\;H_2=\frac{1}{2}
((\frac{a_{41} x_1 +a_{44}x_4}{p_{14}})^2+ ( \frac{ a_{52} x_2+a_{53} x_3}{p_{23}})^2 ).
\end{equation}

{\bf Lie group ${\mathbf A_{4,12}}$:}\\
Inserting $P$ and $P^\prime$ in $(\ref{AA9})$ one can obtain the compatible Poisson structures
${\bf P}$ and ${\bf P^\prime}$ on the Lie group ${\mathbf A_{4,12}}$ as follows:
\begin{equation}
{\bf P}=\left(
\begin{matrix}
 0&0& -p_{14} e^{x_3} sin(x_4) &p_{14} e^{x_3} cos(x_4)\cr
\ast &0&- p_{14} e^{x_3} cos(x_4) &-p_{14} e^{x_3} sin(x_4) \cr
\ast &\ast &0&0\cr
\ast &\ast &\ast &0
 \end{matrix}
\right),
 \end{equation}
 \begin{equation}
  {\bf P^\prime} =
 \left(
\begin{matrix}
0&0&e^{x_3}(-\alpha cos(x_4)+\beta sin(x_4)) & -e^{x_3}(\beta cos(x_4)+\alpha sin(x_4))\cr
\ast &0 &e^{x_3}(\beta cos(x_4)+\alpha sin(x_4)) & e^{x_3}(-\alpha cos(x_4)+\beta sin(x_4)) \cr
\ast &\ast &0&0\cr
\ast &\ast &\ast & 0
 \end{matrix}
\right),
\end{equation}
where $ \alpha =a_{54} x_3-a_{64}x_4$ and $\beta =a_{64} x_3+a_{54}x_4$.
Now by means of $(\ref{A3130})$  and (\ref{A31}), the integrals of motion can be found for this Lie group as follows:
\begin{equation}
H_1=\frac{-2a_{64} x_3 -2a_{54} x_4}{p_{14}} \;,\;\;H_2=\frac{(a_{64}^2- a_{54}^2) x_3^2 + 4a_{54} a_{64} x_3 x_4
+(a_{54}^2- a_{64}^2) x_4^2 }{p_{14}^2}.
\end{equation}

{\bf Lie group ${\mathbf{ A_2 \oplus A_2}}$:}\\
Inserting $P$ and $P^\prime$ in $(\ref{AA9})$ one can obtain the compatible Poisson structures
${\bf P}$ and ${\bf P^\prime}$ on the Lie group ${\mathbf A_2 \oplus A_2}$ as follows:
\begin{equation}
{\bf P}=\left(
\begin{matrix}
 0& p_{12} &0&0 \cr
\ast &0&0& p_{24} \cr
\ast &\ast &0&p_{34}\cr
\ast &\ast &\ast &0
 \end{matrix}
\right),
 \end{equation}
 \begin{equation}
  {\bf P^\prime} =
 \left(
\begin{matrix}
0&a_{21}x_1&0& 0\cr
\ast &0 &0 & p^\prime _{24}+a_{61} x_1+a_{63}x_3 \cr
\ast &\ast &0& a_{73} x_3\cr
\ast &\ast &\ast & 0
 \end{matrix}
\right).
\end{equation}
Now by means of $(\ref{A3130})$  and (\ref{A31}), the  integrals of motion can be found for this Lie group as follows:
\begin{equation}
H_1=\frac{a_{21} x_1}{p_{12}}+\frac{a_{73} x_3}{p_{34}} \;,\;\;H_2=\frac{a_{21}^2 x_1^2}{2 p_{12}^2} +\frac{a_{73}^2 x_3^2}{2 p_{34}^2}.
\end{equation}

{\bf Lie group ${\mathbf{ VII_0\oplus R}}$:}\\
Inserting $P$ and $P^\prime$ in $(\ref{AA9})$ one can obtain the compatible Poisson structures
${\bf P}$ and ${\bf P^\prime}$ on the Lie group ${\mathbf{ VII_0\oplus R}}$ as follows:
\begin{equation}
{\bf P}=\left(
\begin{matrix}
 0& p_{12} &0&0 \cr
\ast &0&0&0 \cr
\ast &\ast &0&p_{34}\cr
\ast &\ast &\ast &0
 \end{matrix}
\right),
 \end{equation}
 \begin{equation}
  {\bf P^\prime} =
 \left(
\begin{matrix}
0&a_{21}x_1+a_{22} x_2&0& 0\cr
\ast &0 &0 & 0 \cr
\ast &\ast &0&a_{73} x_3+a_{74}x_4\cr
\ast &\ast &\ast & 0
 \end{matrix}
\right).
\end{equation}
Now by means of $(\ref{A3130})$  and (\ref{A31}), the integrals of motion can be found for this Lie group as follows:
\begin{equation}
H_1=\frac{a_{21} x_1 +a_{22}x_2}{p_{12}}+\frac{a_{73} x_3+a_{74} x_4}{p_{34}} \;,\;\;H_2=\frac{1}{2}((\frac{a_{21} x_1 +a_{22}x_2}{p_{12}})^2+(\frac{a_{73} x_3+a_{74} x_4}{p_{34}})^2  ).
\end{equation}

{\bf Lie group ${\mathbf{ VI_0\oplus R}}$:}\\
Inserting $P$ and $P^\prime$ in $(\ref{AA9})$ one can obtain the compatible Poisson structures
${\bf P}$ and ${\bf P^\prime}$ on the Lie group ${\mathbf{ VI_0\oplus R}}$ as follows:
\begin{equation}
{\bf P}=\left(
\begin{matrix}
 0& p_{12} &0&0 \cr
\ast &0&0&0 \cr
\ast &\ast &0&p_{34}\cr
\ast &\ast &\ast &0
 \end{matrix}
\right),
 \end{equation}
 \begin{equation}
  {\bf P^\prime} =
 \left(
\begin{matrix}
0&p^\prime _{12}+a_{22}(x_1+ x_2)&0& 0\cr
\ast &0 &0 & 0 \cr
\ast &\ast &0&p^\prime _{34}+a_{73}x_3+a_{74} x_4)\cr
\ast &\ast &\ast & 0
 \end{matrix}
\right).
\end{equation}
Now by means of $(\ref{A3130})$  and (\ref{A31}), the integrals of motion can be found for this Lie group as follows:
\begin{equation}
H_1=\frac{p^\prime _{12}+a_{22} x_1 +a_{22}x_2}{p_{12}}+\frac{p^\prime _{34}+a_{73} x_3+a_{74} x_4}{p_{34}} , H_2=\frac{1}{2}((\frac{p^\prime _{12}+a_{22} x_1 +a_{22}x_2}{p_{12}})^2+(\frac{p^\prime _{34}+a_{73} x_3+a_{74} x_4}{p_{34}})^2  ).
\end{equation}

{\bf Lie group ${\mathbf{ III\oplus R}}$:}\\
Inserting $P$ and $P^\prime$ in $(\ref{AA9})$ one can obtain the compatible Poisson structures
${\bf P}$ and ${\bf P^\prime}$ on the Lie group ${\mathbf{ III\oplus R}}$ as follows:
\begin{equation}
{\bf P}=\left(
\begin{matrix}
 0& p_{13} &p_{13}&0 \cr
\ast &0&0&p_{24} \cr
\ast &\ast &0&0\cr
\ast &\ast &\ast &0
 \end{matrix}
\right),
 \end{equation}
 \begin{equation}
  {\bf P^\prime} =
 \left(
\begin{matrix}
0&a_{33} x_3&a_{33} x_3& 0\cr
\ast &0 &0 & a_{64} x_4\cr
\ast &\ast &0&0\cr
\ast &\ast &\ast & 0
 \end{matrix}
\right).
\end{equation}
Now by means of $(\ref{A3130})$  and (\ref{A31}), the integrals of motion can be found for this Lie group as follows:
\begin{equation}
H_1=\frac{a_{33} x_3 }{p_{13}}+\frac{a_{64} x_4}{p_{24}} \;,\;\;H_2=\frac{1}{2}((\frac{a_{33} x_3 }{p_{13}})^2+(\frac{a_{64} x_4}{p_{24}})^2 ).
\end{equation}

{\bf Lie group ${\mathbf{ II\oplus R}}$:}\\
Inserting $P$ and $P^\prime$ in $(\ref{AA9})$ one can obtain the compatible Poisson structures
${\bf P}$ and ${\bf P^\prime}$ on the Lie group ${\mathbf{ II\oplus R}}$ as follows:
\begin{equation}
{\bf P}=\left(
\begin{matrix}
 0& p_{12} &p_{13}&0 \cr
\ast &0&0&0 \cr
\ast &\ast &0&p_{34}\cr
\ast &\ast &\ast &0
 \end{matrix}
\right),
 \end{equation}
 \begin{equation}
  {\bf P^\prime} =
 \left(
\begin{matrix}
0&a_{44} x_2& p^\prime _{13}+a_{32} x_2-a_{44} x_3& a_{44} x_4\cr
\ast &0 &p^\prime _{23}& 0\cr
\ast &\ast &0&0\cr
\ast &\ast &\ast & 0
 \end{matrix}
\right).
\end{equation}
Now by means of $(\ref{A3130})$  and (\ref{A31}), the integrals of motion can be found for this Lie group as follows:
\begin{equation}
H_1=\frac{a_{44} x_2 }{p_{12}} \;,\;\;H_2=\frac{a_{44}(a_{44} p_{34} x_2^2 - 2 p_{12} p^\prime _{23} x_4)}{2 p_{12}^2 p_{34}}.
\end{equation}

\section {\large {\bf Some compatible Poisson structures on symplectic nilpotent  six dimensional real Lie algebras}}
In this section, we will solve matrix equations (\ref{AA10}), (\ref{AA12})  and (\ref{AA13}) for symplectic nilpotent  six dimensional
 real Lie groups in order to obtain $P$ and $P^\prime $. The list of symplectic real six
dimensional nilpotent Lie algebras is given in appendix B. Note that for calculating $P$ and $P^\prime $ from (\ref{AA12})
 and (\ref{AA13}) we must first calculate vielbeins $e_\alpha ^{\;\mu } $ for nilpotent  6-dimensional  real Lie groups. To this end,
we use the following relation:
\begin{equation}
g^{-1} dg= e_\alpha ^{\;\mu } X_{\mu} dx^{\alpha } \;\;\;\;\;\;g\in {\bf G}
\end{equation}
With the following parameterizations for the real 6-dimensional Lie groups {\bf G}:
\begin{equation}
g = e^{x_1 X_1}e^{x_2 X_2} e^{x_3 X_3} e^{x_4 X_4}e^{x_5 X_5}e^{x_6 X_6}.
\end{equation}
where ${X_i}$ and ${x_i}$ are generators and coordinates of Lie group, respectively. Then,  for left  invariant Lie algebra valued one forms, we have:

\begin{equation}
 g^{-1} dg= dx_1 e^{-x_6 X_6}e^{-x_5 X_5} e^{-x_4 X_4}e^{-x_3 X_3}( e^{-x_2 X_2}X_1 e^{x_2 X_2})e^{x_3 X_3}
 e^{x_4 X_4}e^{x_5 X_5}e^{x_6 X_6}
 \end{equation}
$\hspace*{4cm} +dx_2 e^{-x_6 X_6}e^{-x_5 X_5} e^{-x_4 X_4}( e^{-x_3 X_3}X_2 e^{x_3 X_3})
 e^{x_4 X_4}e^{x_5 X_5}e^{x_6 X_6}$\\ $\hspace*{4cm} +dx_3 e^{-x_6 X_6}e^{-x_5 X_5} ( e^{-x_4 X_4}X_3 e^{x_4 X_4})
 e^{x_5 X_5}e^{x_6 X_6}$\\
$ \hspace*{4cm}+dx_4 e^{-x_6 X_6}( e^{-x_5 X_5}X_4 e^{x_5 X_5})e^{x_6 X_6}+dx_5 e^{-x_6 X_6}X_5 e^{x_6 X_6}
+dx_6 X_6$\\

 such that, for this calculation one can use the following relation \cite{JR}
\begin{equation}
(e^{-x_i X_i}X_j e^{x_i X_i}) = (e^{x_i {\cal X}_i})_j ^{\;k} X_k,
\end{equation}
in which  we have a summation over the index $k$ on the right hand side but there is not summation over the index $i$. In this way one can
calculate all left invariant one forms and vielbeins. 

Let us consider an example for calculating of $P$ and $P^\prime $; for Lie algebra $A_{6,1}$ we have the
following non zero commutators and the matrices ${\cal X}_i$ and ${\cal
Y}^i$:
\begin{equation}
[{\bf e}_1,{\bf e}_2]={\bf e}_3\;,\;\;[{\bf e}_1,{\bf e}_3]={\bf e}_4\;,\;\;[{\bf e}_1,{\bf e}_5]={\bf e}_6\;,
\end{equation}

{\tiny \begin{eqnarray}
{\cal X}_1=\left(
\begin{matrix}
0&0&0&0&0&0\cr
0&0&-1&0&0&0\cr
0&0&0&-1&0&0\cr
0&0&0&0&0&0\cr
0&0&0&0&0&-1\cr
0&0&0&0&0&0
 \end{matrix}
\right) ,\; {\cal X}_2=\left(
\begin{matrix}
0&0&1&0&0&0\cr
0&0&0&0&0&0\cr
0&0&0&0&0&0\cr
0&0&0&0&0&0\cr
0&0&0&0&0&0\cr
0&0&0&0&0&0
 \end{matrix}
\right) ,\; {\cal X}_3=\left(
\begin{matrix}
0&0&0&1&0&0\cr
0&0&0&0&0&0\cr
0&0&0&0&0&0\cr
0&0&0&0&0&0\cr
0&0&0&0&0&0\cr
0&0&0&0&0&0
 \end{matrix}
\right) ,
\nonumber\\
 {\cal X}_4=\left(
\begin{matrix}
0&0&0&0&0&0\cr
0&0&0&0&0&0\cr
0&0&0&0&0&0\cr
0&0&0&0&0&0\cr
0&0&0&0&0&0\cr
0&0&0&0&0&0
 \end{matrix}
\right),\; {\cal X}_5=\left(
\begin{matrix}
0&0&0&0&0&1\cr
0&0&0&0&0&0\cr
0&0&0&0&0&0\cr
0&0&0&0&0&0\cr
0&0&0&0&0&0\cr
0&0&0&0&0&0
 \end{matrix}
\right),\; {\cal X}_6=\left(
\begin{matrix}
0&0&0&0&0&0\cr
0&0&0&0&0&0\cr
0&0&0&0&0&0\cr
0&0&0&0&0&0\cr
0&0&0&0&0&0\cr
0&0&0&0&0&0
 \end{matrix}
\right),\nonumber\\
{\cal Y}^1=\left(
\begin{matrix}
0&0&0&0&0&0\cr
0&0&0&0&0&0\cr
0&0&0&0&0&0\cr
0&0&0&0&0&0\cr
0&0&0&0&0&0\cr
0&0&0&0&0&0
 \end{matrix}
\right) ,\; {\cal Y}^2=\left(
\begin{matrix}
0&0&0&0&0&0\cr
0&0&0&0&0&0\cr
0&0&0&0&0&0\cr
0&0&0&0&0&0\cr
0&0&0&0&0&0\cr
0&0&0&0&0&0
 \end{matrix}
\right) ,\; {\cal Y}^3=\left(
\begin{matrix}
0&-1&0&0&0&0\cr
1&0&0&0&0&0\cr
0&0&0&0&0&0\cr
0&0&0&0&0&0\cr
0&0&0&0&0&0\cr
0&0&0&0&0&0
 \end{matrix}
\right) ,\nonumber\\
 {\cal Y}^4=\left(
\begin{matrix}
0&0&-1&0&0&0\cr
0&0&0&0&0&0\cr
1&0&0&0&0&0\cr
0&0&0&0&0&0\cr
0&0&0&0&0&0\cr
0&0&0&0&0&0
 \end{matrix}
\right),\; {\cal Y}^5=\left(
\begin{matrix}
0&0&0&0&0&0\cr
0&0&0&0&0&0\cr
0&0&0&0&0&0\cr
0&0&0&0&0&0\cr
0&0&0&0&0&0\cr
0&0&0&0&0&0
 \end{matrix}
\right),\; {\cal Y}^6=\left(
\begin{matrix}
0&0&0&0&-1&0\cr
0&0&0&0&0&0\cr
0&0&0&0&0&0\cr
0&0&0&0&0&0\cr
1&0&0&0&0&0\cr
0&0&0&0&0&0
 \end{matrix}
\right).
\end{eqnarray}}

For the Lie algebra $A_{6,1}$ the vielbein matrix $e_\alpha ^{\;\mu} $ has the following form:
\begin{equation}
(e_\alpha ^{\;\mu})=
\left(
\begin{matrix}
1&0&x_2&x_3&0&x_5\cr
0&1&0&0&0&0\cr
0&0&1&0&0&0\cr
0&0&0&1&0&0\cr
0&0&0&0&1&0\cr
0&0&0&0&0&1
 \end{matrix}
\right).
\end{equation}
Now, substituting ${\cal X}_i$, ${\cal Y}^i$ and $e_\alpha ^{\;\mu} $  in (\ref{AA10}), (\ref{AA12})  and (\ref{AA13})  one can obtain
the compatible Poisson structures $P$ and $P^\prime$ for Lie algebra $A_{6,1}$. One of the solutions has the following forms:
\begin{equation}
P=\left(
\begin{matrix}
0&0&0&p _{14} &0&0 \cr
\ast &0&p _{23}&0&0&0\cr
\ast &\ast &0&p_{34}&0&p_{36}\cr
\ast &\ast &\ast &0&0&0\cr
\ast &\ast &\ast &\ast &0&p_{56}\cr
\ast &\ast &\ast &\ast &\ast &0\cr
 \end{matrix}
\right),\;P^\prime=\left(
\begin{matrix}
0&0&0&\frac{c_{44} p _{14}x_4}{p_{34}} &0&0 \cr
\ast &0&b _{32}x_2&0&0&0\cr
\ast &\ast &0&c_{44}x_4&0&0\cr
\ast &\ast &\ast &0&0&0\cr
\ast &\ast &\ast &\ast &0&e_{65} x_5\cr
\ast &\ast &\ast &\ast &\ast &0\cr
 \end{matrix}
\right).
\end{equation}

In this way we have obtained some compatible Poisson structures ($P$ is constant and $P^{\prime}$ as linear function of Lie group
coordinates) on symplectic   nilpotent  six dimensional real Lie algebras, the results are summarized in Table 2.\\

{\small {\bf Table 2}}: {\small  Compatible Poisson structures on nilpotent  6-dimensional real Lie algebra.}\\
    \begin{tabular}{l l l l p{0.15mm} }
    \hline\hline
{\footnotesize ${\bf g}$ }& {\footnotesize   }
&{\footnotesize Non-zero Poisson structure relations }& {\footnotesize Comments}\\ \hline

{\footnotesize $A_{6,1} $}& {\footnotesize $P$} &{\footnotesize $ \lbrace x_1,x_4\rbrace=p_{14}\;,\;\;
\lbrace x_2,x_3\rbrace=p_{23}\;,\;\;\lbrace x_3,x_4\rbrace=p_{34}\;,\;\;\lbrace x_3,x_4\rbrace=p_{36}\;,\;\;
\lbrace x_5,x_6\rbrace=p_{56}$}
&{\footnotesize $p _{14} p _{26}\neq 0$} \\
\vspace{0.7mm}
{\footnotesize }& {\footnotesize $P^\prime $} &{\footnotesize $ \lbrace x_1,x_4\rbrace=\frac{c_{44} p_{14} x_4}{p_{34}}\;
,\;\;\lbrace x_2,x_3\rbrace=b_{32} x_2\;,\;\;\lbrace x_3,x_4\rbrace=c_{44} x_4\;,\;\;\lbrace x_5,x_6\rbrace=e_{65} x_5$}
&{\footnotesize } \\
\vspace{0.7mm}
{\footnotesize $A_{6,7} $}& {\footnotesize $P$} &{\footnotesize $ \lbrace x_1,x_5\rbrace=p_{15}\;,\;\;
\lbrace x_2,x_6\rbrace=p_{26}\;,\;\;\lbrace x_3,x_4\rbrace=p_{34}\;,\;\;\lbrace x_4,x_5\rbrace=p_{45}$}
&{\footnotesize $p _{15} p _{26} p_{34}\neq 0$} \\
\vspace{0.7mm}
{\footnotesize }& {\footnotesize $P^\prime $} &{\footnotesize $ \lbrace x_1,x_5\rbrace=\frac{a_{51}  x_1}{p_{34}}\;
,\;\;\lbrace x_2,x_6\rbrace=b_{62} x_2+b_{66} x_6\;,\;\;\lbrace x_3,x_4\rbrace=c_{43} x_3$}
&{\footnotesize } \\
\vspace*{0.7mm}
{\footnotesize $A_{6,9} $}& {\footnotesize $P$} &{\footnotesize $ \lbrace x_1,x_4\rbrace=p_{14}\;,\;\;
\lbrace x_2,x_6\rbrace=-p_{35}\;,\;\;\lbrace x_3,x_5\rbrace=p_{35}\;,\;\;\lbrace x_3,x_6\rbrace=p_{36}
\;,\;\;\lbrace x_4,x_6\rbrace=p_{46},$}
&{\footnotesize $p _{14} p _{35} \neq 0$} \\
{\footnotesize }& {\footnotesize } &{\footnotesize $ \;\lbrace x_5,x_6\rbrace=p_{56}$}
&{\footnotesize } \\
\vspace{0.7mm}
{\footnotesize }& {\footnotesize $P^\prime $} &{\footnotesize $ \lbrace x_1,x_4\rbrace=a_{41} x_1+
a_{44} x_4+\frac{a_{44}p_{46}  x_2}{p_{35}}\;
\;,\;\;\lbrace x_2,x_6\rbrace=-\frac{d_{62} p_{35}x_2}{p_{46}},$}
&{\footnotesize } \\
{\footnotesize }& {\footnotesize } &{\footnotesize $ \lbrace x_3,x_5\rbrace=\frac{d_{62} p_{35}x_2}{p_{46}}
\;,\;\;\lbrace x_3,x_6\rbrace=c_{62}x_2+c_{65}x_5\;,\;\lbrace x_4,x_6\rbrace=d_{62} x_2\;,\;\lbrace x_5,x_6\rbrace=e_{62}
x_2-\frac{d_{62} p_{35}x_5}{p_{46}}$}
&{\footnotesize } \\
\vspace*{0.7mm}
{\footnotesize $A_{6,24} $}& {\footnotesize $P$} &{\footnotesize $ \lbrace x_1,x_3\rbrace=p_{13}\;,\;\;
\lbrace x_2,x_5\rbrace=p_{25}\;,\;\;\lbrace x_4,x_6\rbrace=p_{46}$}
&{\footnotesize $p _{13} p _{25} p_{46}\neq 0$} \\
\vspace{0.7mm}
{\footnotesize }& {\footnotesize $P^\prime $} &{\footnotesize $ \lbrace x_1,x_3\rbrace=a_{31}  x_1+ a_{33} x_3\;
,\;\;\lbrace x_2,x_5\rbrace=b_{52} x_2+b_{55} x_5\;,\;\;\lbrace x_4,x_6\rbrace=d_{64} x_4+d_{66}x_6$}
&{\footnotesize } \\
\vspace*{0.7mm}
{\footnotesize $A_{6,25} $}& {\footnotesize $P$} &{\footnotesize $ \lbrace x_1,x_3\rbrace=p_{13}\;,\;\;
\lbrace x_2,x_4\rbrace=p_{24}\;,\;\;\lbrace x_3,x_6\rbrace=p_{36}\;,\;\;\lbrace x_5,x_6\rbrace=p_{56}$}
&{\footnotesize $p _{13} p _{24} p_{56}\neq 0$} \\
\vspace{0.7mm}
{\footnotesize }& {\footnotesize $P^\prime $} &{\footnotesize $ \lbrace x_1,x_3\rbrace=a_{31}  x_1\;
,\;\;\lbrace x_2,x_4\rbrace=b_{42} x_2+b_{44} x_4\;,\;\;\lbrace x_5,x_6\rbrace=e_{65} x_5$}
&{\footnotesize } \\
\vspace*{0.7mm}
{\footnotesize $A_{6,26} $}& {\footnotesize $P$} &{\footnotesize $ \lbrace x_1,x_4\rbrace=p_{14}\;,\;\;
\lbrace x_2,x_3\rbrace=p_{23}\;,\;\;\lbrace x_3,x_4\rbrace=\frac{c_{44} p_{14}}{a_{44}}\;,\;\;\lbrace x_3,x_5\rbrace=p_{35}
\;,\;\;\lbrace x_3,x_6\rbrace=\frac{d_{64} p_{35}}{d_{54}},$}
&{\footnotesize $p _{13} p _{24} p_{56}\neq 0$} \\
{\footnotesize }& {\footnotesize } &{\footnotesize $\lbrace x_4,x_5\rbrace=\frac{d_{54} p_{14}}{a_{44}}
\;,\;\;\lbrace x_4,x_6\rbrace=\frac{d_{64} p_{14}}{a_{44}}\;,\;\;\lbrace x_5,x_6\rbrace=p_{56}$}
&{\footnotesize } \\
\vspace{0.7mm}
{\footnotesize }& {\footnotesize $P^\prime $} &{\footnotesize $ \lbrace x_1,x_4\rbrace=a_{44}  x_4\;
,\;\;\lbrace x_2,x_3\rbrace=b_{32} x_2\;,\;\;\lbrace x_3,x_4\rbrace=c_{44} x_4\;,\;\;\lbrace x_4,x_5\rbrace=d_{54} x_4
\;,$}
&{\footnotesize } \\
{\footnotesize }& {\footnotesize } &{\footnotesize $ \lbrace x_4,x_6\rbrace=d_{64} x_4\;,\;\;\lbrace x_5,x_6\rbrace=e_{65} x_5-\frac{d_{54} e_{65} x_6}{d_{64}}$}
&{\footnotesize } \\
\vspace*{0.7mm}
{\footnotesize $A_{6,27} $}& {\footnotesize $P$} &{\footnotesize $ \lbrace x_1,x_2\rbrace=p_{12}\;,\;\;
\lbrace x_1,x_3\rbrace=p_{13}\;,\;\;\lbrace x_2,x_4\rbrace=p_{24}\;,\;\;\lbrace x_5,x_6\rbrace=p_{56}$}
&{\footnotesize $p _{13} p _{24} p_{56}\neq 0$} \\
\vspace{0.7mm}
{\footnotesize }& {\footnotesize $P^\prime $} &{\footnotesize $ \lbrace x_1,x_3\rbrace=a_{33}  x_3\;
,\;\;\lbrace x_2,x_4\rbrace=b_{44} x_4\;,\;\;\lbrace x_5,x_6\rbrace=e_{65} x_5+e_{66} x_6$}
&{\footnotesize } \\
\vspace*{0.7mm}
{\footnotesize $A_{6,32} $}& {\footnotesize $P$} &{\footnotesize $ \lbrace x_1,x_2\rbrace=p_{12}\;,\;\;
\lbrace x_1,x_4\rbrace=-p_{23}\;,\;\;\lbrace x_1,x_6\rbrace=\frac{a_{63}p_{23}+a_{25}p_{56}}{b_{33}}\;,\;\;\lbrace x_2,x_3\rbrace=p_{23}\;,\;\;\lbrace x_5,x_6\rbrace=p_{56}$}
&{\footnotesize $p _{23} p _{56} \neq 0$} \\
\vspace{0.7mm}
{\footnotesize }& {\footnotesize $P^\prime $} &{\footnotesize $ \lbrace x_1,x_2\rbrace=-b_{33}  x_1-b_{34} x_2+a_{23} x_3
+a_{24} x_4+a_{25} x_5\;,\;\;\lbrace x_1,x_4\rbrace=-b_{33} x_3-b_{34} x_4,$}
&{\footnotesize } \\
{\footnotesize }& {\footnotesize } &{\footnotesize $ \lbrace x_1,x_6\rbrace=a_{63} x_3+\frac{a_{63}(b_{33}b_{34}+a_{25}e_{66})x_4}{b_{33}^2}+\frac{a_{25}e_{65}
x_5+a_{25}e_{66}x_6}{b_{33}}\;,\;\;\lbrace x_2,x_3\rbrace=b_{33} x_3+b_{34}x_4
\;,$}
&{\footnotesize } \\
{\footnotesize }& {\footnotesize } &{\footnotesize $ \lbrace x_5,x_6\rbrace=\frac{a_{63} e_{66}x_4}{b_{33}}+e_{65} x_5+e_{66} x_{6}$}
&{\footnotesize } \\

\hline\hline
\end{tabular}

\section {\large {\bf Integrable bi-Hamiltonian systems on symplectic nilpotent  six dimensional real Lie groups}}

In this section, we construct the integrable bi-Hamiltonian systems
with symplectic nilpotent  six dimensional real Lie groups as phase space.
Using the matrices $P$ and $P^\prime$ given in Table 2 and inserting $P$ and $P^\prime$ in $(\ref{AA9})$
and using the related vielbeins as in appendix B the compatible Poisson structures  ${\bf P}$ and ${\bf P^\prime}$
on Lie groups are obtained, then using  $(\ref{A3130})$ and $(\ref{A31})$ one can find the  Hamiltonian and integrals
of motions of bi-Hamiltonian systems. In the following we   perform this work separately   for symplectic nilpotent  six dimensional real Lie groups. In this way we  obtain new bi-Hamiltonian  systems over nilpotent  six dimensional real Lie groups as phase spaces.

{\bf Lie group ${\mathbf A_{6,1}}$:}\\
Inserting $P$,  $P^\prime$ and vielbeins matrix in $(\ref{AA9})$ one can obtain the compatible Poisson structures
${\bf P}$ and ${\bf P^\prime}$ on the Lie group ${\mathbf A_{6,1}}$ as follows:
\begin{equation}
{\bf P}=\left(
\begin{matrix}
0&0&0&p_{14}&0 & 0 \cr
\ast &0&p_{23} &0&0&0 \cr
\ast &\ast &0&p_{34}+p_{14} x_2&0&p_{36}\cr
\ast &\ast &\ast &0&0&-p_{14} x_5\cr
\ast &\ast &\ast &\ast &0&p_{56} \cr
\ast &\ast &\ast &\ast &\ast &0
 \end{matrix}
\right),
 \end{equation}
 \begin{equation}
  {\bf P^\prime} =\left(
\begin{matrix}
0&0&0&\frac{c_{44}p_{14}x_4}{p_{34}}&0 & 0 \cr
\ast &0&b_{32}x_2 &0&0&0 \cr
\ast &\ast &0&c_{44}x_4+\frac{c_{44} p_{14} x_2 x_4}{p_{34}}&0&0\cr
\ast &\ast &\ast &0&0&-\frac{c_{44} p_{14}x_4 x_5}{p_{34}}\cr
\ast &\ast &\ast &\ast &0&e_{65}x_5 \cr
\ast &\ast &\ast &\ast &\ast &0
 \end{matrix}
\right).
\end{equation}
Now by means of $(\ref{A3130})$  and $(\ref{A31})$, the integrals of motion can be found for this Lie group as follows:
\begin{equation*}
H_1=\frac{b_{32} x_2}{p_{23}}+\frac{c_{44} x_4}{p_{34}}+\frac{e_{65} x_5}{p_{56}} \;,\;\;H_2=\frac{1}{2}((\frac{b_{32} x_2}{p_{23}})^2+(\frac{c_{44} x_4}{p_{34}})^2+(\frac{e_{65} x_5}{p_{56}})^2)\;,
\end{equation*}
\begin{equation}
H_3=\frac{1}{3}((\frac{b_{32} x_2}{p_{23}})^3+(\frac{c_{44} x_4}{p_{34}})^3+(\frac{e_{65} x_5}{p_{56}})^3).
\end{equation}

{\bf Lie group ${\mathbf A_{6,7}}$:}\\
Inserting $P$,  $P^\prime$ and vielbeins matrix in $(\ref{AA9})$ one can obtain the compatible Poisson structures
${\bf P}$ and ${\bf P^\prime}$ on the Lie group ${\mathbf A_{6,7}}$ as follows:
\begin{equation}
{\bf P}=\left(
\begin{matrix}
0&0&0&0&p_{15} & 0 \cr
\ast &0&0 &0&0&p_{26} \cr
\ast &\ast &0&p_{34}&0&0\cr
\ast &\ast &\ast &0&p_{45}+p_{15} x_3&0\cr
\ast &\ast &\ast &\ast &0&0 \cr
\ast &\ast &\ast &\ast &\ast &0
 \end{matrix}
\right),
 \end{equation}
 \begin{equation}
  {\bf P^\prime} = \left(
\begin{matrix}
0&0&0&0&a_{51}x_1& 0\cr
\ast &0&0 &0&0&b_{62}x_2+b_{66} x_6 \cr
\ast &\ast &0&c_{43}x_3&0&0\cr
\ast &\ast &\ast &0&a_{51}x_1 x_3&0\cr
\ast &\ast &\ast &\ast &0&0 \cr
\ast &\ast &\ast &\ast &\ast &0
 \end{matrix}
\right).
\end{equation}
Now by means of (\ref{A31}), the integrals of motion can be found for this Lie group as follows:
\begin{equation*}
H_1=\frac{a_{51} x_1}{p_{15}}+\frac{b_{62} x_2}{p_{26}}+\frac{c_{43} x_3}{p_{34}}+\frac{b_{66} x_6}{p_{26}} \;,\;\;H_2=\frac{1}{2}((\frac{a_{51} x_1}{p_{15}})^2+(\frac{c_{43} x_3}{p_{34}})^2+(\frac{b_{62}x_2+b_{66} x_6}{p_{26}})^2)\;,
\end{equation*}
\begin{equation}
H_3=\frac{1}{3}((\frac{a_{51} x_1}{p_{15}})^3+(\frac{c_{43} x_3}{p_{34}})^3+(\frac{b_{62}x_2+b_{66} x_6}{p_{26}})^3).
\end{equation}

{\bf Lie group ${\mathbf A_{6,9}}$:}\\
Inserting $P$,  $P^\prime$ and vielbeins matrix in $(\ref{AA9})$ one can obtain the compatible Poisson structures
${\bf P}$ and ${\bf P^\prime}$ on the Lie group ${\mathbf A_{6,9}}$ as follows:
\begin{equation}
{\bf P}=\left(
\begin{matrix}
0&0&0&p_{14}&0 & 0 \cr
\ast &0&0 &0&0&-p_{35} \cr
\ast &\ast &0&p_{14}x_2&p_{35}&p_{36}\cr
\ast &\ast &\ast &0&0&p_{46}+1/2 p_{14}x_2^2-p_{14} x_5\cr
\ast &\ast &\ast &\ast &0&p_{56} \cr
\ast &\ast &\ast &\ast &\ast &0
 \end{matrix}
\right),
 \end{equation}
 \begin{equation}
  {\bf P^\prime} = {\tiny \left(
\begin{matrix}
0&0&0&a_{41} x_1+\frac{a_{44} p_{46} x_2}{p_{35}}+a_{44} x_4&0& 0\cr
\ast &0&0 &0&0&-\frac{d_{62} p_{35}x_2}{p_{46}} \cr
\ast &\ast &0&a_{41}x_1 x_2+\frac{a_{44}p_{46}x_2^2}{p_{35}}+a_{44} x_2 x_4&\frac{d_{62} p_{35}x_2}{p_{46}}&
c_{62}x_2+c_{65}x_5\cr
\ast &\ast &\ast &0&0&\frac{2d_{62} p_{35}x_2+(a_{41}p_{35}x_1+a_{44}p_{46}x_2+a_{44}p_{35}x_4)(x_2^2-2x_5)}{2p_{35}}\cr
\ast &\ast &\ast &\ast &0&e_{62} x_2-\frac{d_{62}p_{35}x_5}{p_{46}} \cr
\ast &\ast &\ast &\ast &\ast &0
 \end{matrix}
\right)}.
\end{equation}
Now by means of (\ref{A31}), the integrals of motion can be found for this Lie group as follows:
\begin{equation*}
H_1=\frac{a_{41}p_{35}p_{46}x_1+2d_{62}p_{14}p_{35}x_2+a_{44}p_{46}(p_{46} x_2+p_{35}x_4)}{p_{14}p_{35}p_{46}} \;,\;\;
\end{equation*}
\begin{equation*}
H_2=\frac{1}{4}((\frac{4 d_{62}x_2}{p_{46}})^2+(\frac{a_{41}p_{35}x_1+a_{44}p_{46}x_2+a_{44}p_{35}x_4}{p_{14}p_{35}} )^2+(\frac{a_{41}x_1+a_{44}(\frac{p_{46}x_2}{p_{35}}+x_4) }{p_{14}} )^2)\;,
\end{equation*}
\begin{equation}
H_3=\frac{1}{6}((\frac{4 d_{62}x_2}{p_{46}})^3+(\frac{a_{41}p_{35}x_1+a_{44}p_{46}x_2+a_{44}p_{35}x_4}{p_{14}p_{35}} )^3+(\frac{a_{41}x_1+a_{44}(\frac{p_{46}x_2}{p_{35}}+x_4) }{p_{14}} )^3).
\end{equation}

{\bf Lie group ${\mathbf A_{6,24}}$:}\\
Inserting $P$,  $P^\prime$ and vielbeins matrix in $(\ref{AA9})$ one can obtain the compatible Poisson structures
${\bf P}$ and ${\bf P^\prime}$ on the Lie group ${\mathbf A_{6,24}}$ as follows:
\begin{equation}
{\bf P}=\left(
\begin{matrix}
0&0&&p_{13}&0&0 & 0 \cr
\ast &0&0 &0&p_{25}&0 \cr
\ast &\ast &0&0&0&0\cr
\ast &\ast &\ast &0&0&p_{46}\cr
\ast &\ast &\ast &\ast &0&0 \cr
\ast &\ast &\ast &\ast &\ast &0
 \end{matrix}
\right),
 \end{equation}
 \begin{equation}
  {\bf P^\prime} = \left(
\begin{matrix}
0&0&a_{31}x_1+a_{33}x_3&0&0& 0\cr
\ast &0&0 &0&b_{52}x_2+b_{55}x_5&0 \cr
\ast &\ast &0&0&0&0\cr
\ast &\ast &\ast &0&0&d_{64}x_4+d_{66}x_6\cr
\ast &\ast &\ast &\ast &0&0 \cr
\ast &\ast &\ast &\ast &\ast &0
 \end{matrix}
\right).
\end{equation}
Now by means of (\ref{A31}), the integrals of motion can be found for this Lie group as follows:
\begin{equation*}
H_1=\frac{a_{31} x_1+a_{33} x_3}{p_{13}}+\frac{b_{52} x_2+b_{55}x_5}{p_{25}}+\frac{d_{64} x_4+d_{66}x_6}{p_{45}}\;,
\end{equation*}
\begin{equation*}
H_2=\frac{1}{2}((\frac{a_{31} x_1+a_{33} x_3}{p_{13}})^2+(\frac{b_{52} x_2+b_{55}x_5}{p_{25}})^2+(\frac{d_{64} x_4+d_{66}x_6}{p_{45}})^2)\;,
\end{equation*}
\begin{equation}
H_3=\frac{1}{3}((\frac{a_{31} x_1+a_{33} x_3}{p_{13}})^3+(\frac{b_{52} x_2+b_{55}x_5}{p_{25}})^3+(\frac{d_{64} x_4+d_{66}x_6}{p_{45}})^3).
\end{equation}

{\bf Lie group ${\mathbf A_{6,25}}$:}\\
Inserting $P$,  $P^\prime$ and vielbeins matrix in $(\ref{AA9})$ one can obtain the compatible Poisson structures
${\bf P}$ and ${\bf P^\prime}$ on the Lie group ${\mathbf A_{6,25}}$ as follows:
\begin{equation}
{\bf P}=\left(
\begin{matrix}
0&0&p_{13}&0&0 & 0 \cr
\ast &0&0 &p_{24}&0&p_{24}x_5 \cr
\ast &\ast &0&0&0&p_{36}\cr
\ast &\ast &\ast &0&0&0\cr
\ast &\ast &\ast &\ast &0&p_{56} \cr
\ast &\ast &\ast &\ast &\ast &0
 \end{matrix}
\right),
 \end{equation}
 \begin{equation}
  {\bf P^\prime} = \left(
\begin{matrix}
0&0&a_{31}x_1&0&0& 0\cr
\ast &0&0 &b_{42}x_2+b_{44}x_4&0&b_{42}x_2x_5+b_{44}x_4x_5 \cr
\ast &\ast &0&0&0&0\cr
\ast &\ast &\ast &0&0&0\cr
\ast &\ast &\ast &\ast &0&e_{65}x_5\cr
\ast &\ast &\ast &\ast &\ast &0
 \end{matrix}
\right).
\end{equation}
Now by means of (\ref{A31}), the integrals of motion can be found for this Lie group as follows:
\begin{equation*}
H_1=\frac{a_{31} x_1}{p_{13}}+\frac{b_{42} x_2+b_{44}x_4}{p_{24}}+\frac{e_{65} x_5}{p_{56}}\;,\;\;H_2=\frac{1}{2}((\frac{a_{31} x_1}{p_{13}})^2+(\frac{b_{42} x_2+b_{44}x_4}{p_{24}})^2+(\frac{e_{65} x_5}{p_{56}})^2)\;,
\end{equation*}
\begin{equation}
H_3=\frac{1}{3}((\frac{a_{31} x_1}{p_{13}})^3+(\frac{b_{42} x_2+b_{44}x_4}{p_{24}})^3+(\frac{e_{65} x_5}{p_{56}})^3).
\end{equation}

{\bf Lie group ${\mathbf A_{6,26}}$:}\\
Inserting $P$,  $P^\prime$ and vielbeins matrix in $(\ref{AA9})$ one can obtain the compatible Poisson structures
${\bf P}$ and ${\bf P^\prime}$ on the Lie group ${\mathbf A_{6,26}}$ as follows:
\begin{equation}
{\bf P}=\left(
\begin{matrix}
0&0&0&p_{14}&0&0  \cr
\ast &0&p_{23} &0&0&0\cr
\ast &\ast &0&\frac{c_{44}p_{14}}{a_{44}}+p_{14}x_2&p_{35}&\frac{d_{64}p_{35}}{d_{54}}\cr
\ast &\ast &\ast &0&\frac{d_{54}p_{14}}{a_{44}}&\frac{d_{64}p_{14}}{a_{44}}\cr
\ast &\ast &\ast &\ast &0&p_{56} \cr
\ast &\ast &\ast &\ast &\ast &0
 \end{matrix}
\right),
 \end{equation}
 \begin{equation}
  {\bf P^\prime} = \left(
\begin{matrix}
0&0&0&a_{44} x_4&0& 0\cr
\ast &0&b_{32} x_2 &0&0&0 \cr
\ast &\ast &0&c_{44} x_4+a_{44} x_2 x_4&0&0\cr
\ast &\ast &\ast &0&d_{54} x_4&d_{64} x_4\cr
\ast &\ast &\ast &\ast &0&e_{65}x_5-\frac{d_{54}e_{65}x_6}{d_{64}}\cr
\ast &\ast &\ast &\ast &\ast &0
 \end{matrix}
\right).
\end{equation}
Now by means of (\ref{A31}), the integrals of motion can be found for this Lie group as follows:
\begin{equation*}
H_1=\frac{b_{32} x_2}{p_{23}}+\frac{a_{44} x_4}{p_{14}}+\frac{e_{65} (d_{64}x_5-d_{54}x_6)}{d_{64}p_{56}}\;,\;\;H_2=\frac{1}{2}((\frac{b_{32} x_2}{p_{23}})^2+(\frac{a_{44} x_4}{p_{14}})^2+(\frac{e_{65} (d_{64}x_5-d_{54}x_6)}{d_{64}p_{56}})^2)\;,
\end{equation*}
\begin{equation}
H_3=\frac{1}{3}((\frac{b_{32} x_2}{p_{23}})^3+(\frac{a_{44} x_4}{p_{14}})^3+(\frac{e_{65} (d_{64}x_5-d_{54}x_6)}{d_{64}p_{56}})^3).
\end{equation}

{\bf Lie group ${\mathbf A_{6,27}}$:}\\
Inserting $P$,  $P^\prime$ and vielbeins matrix in $(\ref{AA9})$ one can obtain the compatible Poisson structures
${\bf P}$ and ${\bf P^\prime}$ on the Lie group ${\mathbf A_{6,27}}$ as follows:
\begin{equation}
{\bf P}=\left(
\begin{matrix}
0&p_{12}&p_{13}&0&0 & 0 \cr
\ast &0&0&p_{24}&0&0\cr
\ast &\ast &0&0&0&0\cr
\ast &\ast &\ast &0&0&0\cr
\ast &\ast &\ast &\ast &0&p_{56} \cr
\ast &\ast &\ast &\ast &\ast &0
 \end{matrix}
\right),
 \end{equation}
 \begin{equation}
  {\bf P^\prime} = \left(
\begin{matrix}
0&0&a_{33}x_3&0&0& 0\cr
\ast &0&0 &b_{44}x_4&0&0 \cr
\ast &\ast &0&0&0&0\cr
\ast &\ast &\ast &0&0\cr
\ast &\ast &\ast &\ast &0&e_{65}x_5+e_{66}x_6\cr
\ast &\ast &\ast &\ast &\ast &0
 \end{matrix}
\right).
\end{equation}
Now by means of (\ref{A31}), the integrals of motion can be found for this Lie group as follows:
\begin{equation*}
H_1=\frac{a_{33} x_3}{p_{13}}+\frac{b_{44} x_4}{p_{24}}+\frac{e_{65}x_5+e_{66}x_6 }{p_{56}}\;,\;\;H_2=\frac{1}{2}((\frac{a_{33} x_3}{p_{13}})^2+(\frac{b_{44} x_4}{p_{24}})^2+(\frac{e_{65}x_5+e_{66}x_6 }{p_{56}}^2))\;,
\end{equation*}
\begin{equation}
H_3=\frac{1}{3}((\frac{a_{33} x_3}{p_{13}})^3+(\frac{b_{44} x_4}{p_{24}})^3+(\frac{e_{65}x_5+e_{66}x_6 }{p_{56}})^3).
\end{equation}

\section{Concluding remarks}
Using a procedure based on the adjoint representation of the Lie algebra we give a method for
calculation of compatible Poisson structures on four
and nilpotent six dimensional symplectic real Lie algebras. Also by use of Magri-Morosi's theorem we have obtained new
bi-Hamiltonian systems with  these Lie groups as  phase spaces. It seems that using of the adjoint representation of the Lie algebra
(mainly for low dimensions) is simpler than the Lax-pair method (at least in computation). As an open problem, one can obtain another set of  compatible
Poisson structures by setting $P^\prime $ as a second order functions of the Lie group parameters and in this way obtain new bi-Hamiltonian systems.\\\vspace*{.5cm}

{\bf Appendix A:} The list of symplectic four dimensional real Lie algebras \cite{JAF2}\\
    \begin{tabular}{l l | l l }
    \hline\hline
{\footnotesize ${\bf g}$ }&{\footnotesize Non-zero commutation relations }&{\footnotesize ${\bf g}$ }&{\footnotesize Non-zero commutation relations } \\ \hline

{\footnotesize$A_{4,1}$}&{\footnotesize$[{\bf e}_2,{\bf e}_4]={\bf e}_1 \;,\;\;[{\bf e}_3,{\bf e}_4]={\bf e}_2$} &{\footnotesize$A_{4,2}^{-1}$}&{\footnotesize$[{\bf e}_1,{\bf e}_4]=-{\bf e}_1 ,\;[{\bf e}_2,{\bf e}_4]={\bf e}_2 ,\;[{\bf e}_3,{\bf e}_4]={\bf e}_2+{\bf e}_3$} \\

{\footnotesize$A_{4,3}$}&{\footnotesize$[{\bf e}_1,{\bf e}_4]={\bf e}_1 \;,\;\;[{\bf e}_3,{\bf e}_4]={\bf e}_2$} &{\footnotesize$A_{4,5}^{-1,-1}$}&{\footnotesize$[{\bf e}_1,{\bf e}_4]={\bf e}_1 \;,\;\;[{\bf e}_2,{\bf e}_4]=-{\bf e}_2 \;,\;\;[{\bf e}_3,{\bf e}_4]=-{\bf e}_3$} \\

{\footnotesize$A_{4,5}^{-1,b}$}&{\footnotesize$[{\bf e}_1,{\bf e}_4]={\bf e}_1 , [{\bf e}_2,{\bf e}_4]=-{\bf e}_2 , [{\bf e}_3,{\bf e}_4]=b \;{\bf e}_3$} &{\footnotesize$A_{4,5}^{a,-1}$}&{\footnotesize$[{\bf e}_1,{\bf e}_4]={\bf e}_1 \;,\;\;[{\bf e}_2,{\bf e}_4]=a\;{\bf e}_2 \;,\;\;[{\bf e}_3,{\bf e}_4]=-{\bf e}_3$} \\

{\footnotesize$A_{4,5}^{a,-a}$}&{\footnotesize$[{\bf e}_1,{\bf e}_4]={\bf e}_1 , [{\bf e}_2,{\bf e}_4]=a {\bf e}_2 , [{\bf e}_3,{\bf e}_4]=-a {\bf e}_3$} &{\footnotesize$A_{4,6}^{a,0}$}&{\footnotesize$[{\bf e}_1,{\bf e}_4]=a {\bf e}_1 \;,\;\;[{\bf e}_2,{\bf e}_4]=-{\bf e}_3 \;,\;\;[{\bf e}_3,{\bf e}_4]={\bf e}_2$} \\

{\footnotesize$A_{4,7}$}&{\footnotesize$[{\bf e}_1,{\bf e}_4]=2{\bf e}_1 ,\;[{\bf e}_2,{\bf e}_3]={\bf e}_1 ,\;[{\bf e}_2,{\bf e}_4]={\bf e}_2$} &{\footnotesize$A_{4,9}^{0}$}&{\footnotesize$[{\bf e}_1,{\bf e}_4]={\bf e}_1 \;,\;\;[{\bf e}_2,{\bf e}_3]={\bf e}_1 \;,\;\;[{\bf e}_2,{\bf e}_4]={\bf e}_2$} \\

&{\footnotesize$\;[{\bf e}_3,{\bf e}_4]={\bf e}_2+{\bf e}_3$} &{\footnotesize$A_{4,9}^{-1/2}$}&{\footnotesize$[{\bf e}_1,{\bf e}_4]=1/2 {\bf e}_1 \;,\;\;[{\bf e}_2,{\bf e}_3]={\bf e}_1 \;,\;\;[{\bf e}_2,{\bf e}_4]={\bf e}_2$} \\

{\footnotesize$A_{4,9}^1$}&{\footnotesize$[{\bf e}_1,{\bf e}_4]=2{\bf e}_1 ,\;[{\bf e}_2,{\bf e}_3]={\bf e}_1 ,\;[{\bf e}_2,{\bf e}_4]={\bf e}_2$} &&{\footnotesize$[{\bf e}_3,{\bf e}_4]=-1/2{\bf e}_3 $} \\

&{\footnotesize$\;[{\bf e}_3,{\bf e}_4]={\bf e}_3$} &{\footnotesize$A_{4,9}^{b}$}&{\footnotesize$[{\bf e}_1,{\bf e}_4]=(1+b) {\bf e}_1 ,\;[{\bf e}_2,{\bf e}_3]={\bf e}_1 ,\;[{\bf e}_2,{\bf e}_4]={\bf e}_2$} \\

{\footnotesize$A_{4,11}^b$}&{\footnotesize$[{\bf e}_1,{\bf e}_4]=2b {\bf e}_1 \;,\;\;[{\bf e}_2,{\bf e}_3]={\bf e}_1 $} &&{\footnotesize$[{\bf e}_3,{\bf e}_4]=b {\bf e}_3 $} \\

&{\footnotesize$\;[{\bf e}_2,{\bf e}_4]=b {\bf e}_2-{\bf e}_3\;,\;\;[{\bf e}_3,{\bf e}_4]={\bf e}_2+b{\bf e}_3$} &{\footnotesize$A_{4,12}$}&{\footnotesize$[{\bf e}_1,{\bf e}_3]= {\bf e}_1 \;,\;\;[{\bf e}_1,{\bf e}_4]=-{\bf e}_2 \;,\;\;[{\bf e}_2,{\bf e}_3]={\bf e}_2$} \\

{\footnotesize$A_2\oplus A_2$}&{\footnotesize$[{\bf e}_1,{\bf e}_2]= {\bf e}_2 \;,\;\;[{\bf e}_3,{\bf e}_4]={\bf e}_4 $} &&{\footnotesize$[{\bf e}_2,{\bf e}_4]= {\bf e}_1 $} \\

{\footnotesize$VI_0\oplus R$}&{\footnotesize$[{\bf e}_1,{\bf e}_3]={\bf e}_2 \;,\;\;[{\bf e}_2,{\bf e}_3]={\bf e}_1 $} &{\footnotesize$III\oplus R$}&{\footnotesize$[{\bf e}_1,{\bf e}_2]=-{\bf e}_2-{\bf e}_3 \;,\;\;[{\bf e}_1,{\bf e}_3]=-{\bf e}_2-{\bf e}_3 $} \\

{\footnotesize$VII_0\oplus R$}&{\footnotesize$[{\bf e}_1,{\bf e}_3]=-{\bf e}_2 \;,\;\;[{\bf e}_2,{\bf e}_3]={\bf e}_1 $} &{\footnotesize$II\oplus R$}&{\footnotesize$[{\bf e}_2,{\bf e}_3]={\bf e}_1 $} \\

\smallskip \\
\hline\hline
\end{tabular}

\vspace*{.75cm}

{\bf Appendix B:} The list of symplectic nilpotent  6-dimensional real Lie algebras and vielbeins matrix\\
    \begin{tabular}{l l l}
    \hline\hline
{\footnotesize ${\bf g}$ }&{\footnotesize Non-zero commutation relations } &{\footnotesize vielbeins matrix }\\ \hline

{\footnotesize $A_{6,1}$}&{\footnotesize $\begin{matrix}[{\bf e}_1,{\bf e}_2]={\bf e}_3 \cr
[{\bf e}_1,{\bf e}_3]={\bf e}_4 \cr
[{\bf e}_1,{\bf e}_5]={\bf e}_6 \end{matrix}$
}& {\footnotesize
$\left(
\begin{matrix}
1&0&x_2&x_3&0& x_5\cr
0 &1&0 &0&0&0 \cr
0 &0 &1&0&0&0\cr
0 &0 &0 &1&0&0\cr
0 &0 &0 &0 &1&0\cr
0 &0 &0 &0&0 &1
 \end{matrix}
\right)$ }\\

{\footnotesize $A_{6,2}$}&{\footnotesize $\begin{matrix}[{\bf e}_1,{\bf e}_2]={\bf e}_3 \cr [{\bf e}_1,{\bf e}_3]={\bf e}_4\cr
[{\bf e}_1,{\bf e}_4]={\bf e}_5 \cr
[{\bf e}_1,{\bf e}_5]={\bf e}_6 \end{matrix}$} & {\footnotesize
$\left(
\begin{matrix}
1&0&x_2&x_3&x_4& x_5\cr
0 &1&0 &0&0&0 \cr
0 &0 &1&0&0&0\cr
0 &0 &0 &1&0&0\cr
0 &0 &0 &0 &1&0\cr
0 &0 &0 &0&0 &1
 \end{matrix}
\right)$ }\\

{\footnotesize $A_{6,3}$}&{\footnotesize $\begin{matrix}[{\bf e}_1,{\bf e}_2]={\bf e}_6 \cr [{\bf e}_1,{\bf e}_3]={\bf e}_4\cr
[{\bf e}_2,{\bf e}_3]={\bf e}_5 \end{matrix}$} & {\footnotesize
$\left(
\begin{matrix}
1&0&0&x_3&0& x_2\cr
0 &1&0 &0&x_3&0 \cr
0 &0 &1&0&0&0\cr
0 &0 &0 &1&0&0\cr
0 &0 &0 &0 &1&0\cr
0 &0 &0 &0&0 &1
 \end{matrix}
\right)$ }\\

{\footnotesize $A_{6,4}$}&{\footnotesize $\begin{matrix} [{\bf e}_1,{\bf e}_2]={\bf e}_5 \cr [{\bf e}_1,{\bf e}_3]={\bf e}_6\cr
[{\bf e}_2,{\bf e}_4]={\bf e}_6 \end{matrix}$} & {\footnotesize
$\left(
\begin{matrix}
1&0&0&0&x_2& x_3\cr
0 &1&0 &0&0&x_4 \cr
0 &0 &1&0&0&0\cr
0 &0 &0 &1&0&0\cr
0 &0 &0 &0 &1&0\cr
0 &0 &0 &0&0 &1
 \end{matrix}
\right)$ }\\

\smallskip \\
\hline\hline
\end{tabular}

  \begin{tabular}{l l l}
    \hline\hline
{\footnotesize ${\bf g}$ }&{\footnotesize Non-zero commutation relations } &{\footnotesize vielbeins matrix }\\ \hline

{\footnotesize $A_{6,5}$}&{\footnotesize $\begin{matrix} [{\bf e}_1,{\bf e}_3]={\bf e}_5 \;\;\; \cr [{\bf e}_1,{\bf e}_4]={\bf e}_6\;\;\;\cr
[{\bf e}_2,{\bf e}_3]=a\;{\bf e}_6 \cr
[{\bf e}_2,{\bf e}_4]={\bf e}_5\;\;\; \cr  a=\pm 1 \end{matrix} $} &  {\footnotesize
$\left(
\begin{matrix}
1&0&0&0&x_3& x_4\cr
0 &1&0 &0&x_4&a x_3 \cr
0 &0 &1&0&0&0\cr
0 &0 &0 &1&0&0\cr
0 &0 &0 &0 &1&0\cr
0 &0 &0 &0&0 &1
 \end{matrix}
\right)$ }\\

{\footnotesize $A_{6,6}$}&{\footnotesize $\begin{matrix}[{\bf e}_1,{\bf e}_2]={\bf e}_6 \cr [{\bf e}_1,{\bf e}_3]={\bf e}_4\cr
[{\bf e}_1,{\bf e}_4]={\bf e}_5 \cr
[{\bf e}_2,{\bf e}_3]={\bf e}_5 \end{matrix}$} &  {\footnotesize
$\left(
\begin{matrix}
1&0&0&x_3&x_4& x_5\cr
0 &1&0 &0&x_3&0 \cr
0 &0 &1&0&0&0\cr
0 &0 &0 &1&0&0\cr
0 &0 &0 &0 &1&0\cr
0 &0 &0 &0&0 &1
 \end{matrix}
\right)$ }\\

{\footnotesize $A_{6,7}$}&{\footnotesize $\begin{matrix}[{\bf e}_1,{\bf e}_3]={\bf e}_4 \cr [{\bf e}_1,{\bf e}_4]={\bf e}_5\cr
[{\bf e}_2,{\bf e}_3]={\bf e}_6 \end{matrix} $} &  {\footnotesize
$\left(
\begin{matrix}
1&0&0&x_3&x_4& 0\cr
0 &1&0 &0&0&x_3 \cr
0 &0 &1&0&0&0\cr
0 &0 &0 &1&0&0\cr
0 &0 &0 &0 &1&0\cr
0 &0 &0 &0&0 &1
 \end{matrix}
\right)$ }\\

{\footnotesize $A_{6,8}$}&{\footnotesize $\begin{matrix} [{\bf e}_1,{\bf e}_2]={\bf e}_3+{\bf e}_5 \cr
 [{\bf e}_1,{\bf e}_3]={\bf e}_4 \hspace*{0.65cm}\cr
[{\bf e}_2,{\bf e}_5]={\bf e}_6 \hspace*{0.65cm} \end{matrix} $} &  {\footnotesize
$\left(
\begin{matrix}
1&0&x_2&x_3&x_2& -x_2^2/2\cr
0 &1&0 &0&0&x5 \cr
0 &0 &1&0&0&0\cr
0 &0 &0 &1&0&0\cr
0 &0 &0 &0 &1&0\cr
0 &0 &0 &0&0 &1
 \end{matrix}
\right)$ }\\

{\footnotesize $A_{6,9}$}&{\footnotesize $\begin{matrix} [{\bf e}_1,{\bf e}_2]={\bf e}_3 \cr [{\bf e}_1,{\bf e}_3]={\bf e}_4\cr
[{\bf e}_1,{\bf e}_5]={\bf e}_6 \cr
[{\bf e}_2,{\bf e}_3]={\bf e}_6 \end{matrix}$} &  {\footnotesize
$\left(
\begin{matrix}
1&0&x_2&x_3&0&-x_2^2/2+ x_5\cr
0 &1&0 &0&0&x_3 \cr
0 &0 &1&0&0&0\cr
0 &0 &0 &1&0&0\cr
0 &0 &0 &0 &1&0\cr
0 &0 &0 &0&0 &1
 \end{matrix}
\right)$ }\\

{\footnotesize $A_{6,10}$}&{\footnotesize $\begin{matrix} [{\bf e}_1,{\bf e}_2]={\bf e}_3 \;\;\;\cr
 [{\bf e}_1,{\bf e}_3]={\bf e}_5 \;\;\;\cr
[{\bf e}_1,{\bf e}_4]={\bf e}_6\;\;\;\cr
[{\bf e}_2,{\bf e}_3]= a\;{\bf e}_6 \cr
[{\bf e}_2,{\bf e}_4]= {\bf e}_5 \;\;\;\cr a=\pm 1 \end{matrix}$} &  {\footnotesize
$\left(
\begin{matrix}
1&0&x_2&0&x_3& -a x_2^2/2+x_4\cr
0 &1&0 &0&x_4&a x_3 \cr
0 &0 &1&0&0&0\cr
0 &0 &0 &1&0&0\cr
0 &0 &0 &0 &1&0\cr
0 &0 &0 &0&0 &1
 \end{matrix}
\right)$ }\\

{\footnotesize $A_{6,11}$}&{\footnotesize $\begin{matrix} [{\bf e}_1,{\bf e}_2]={\bf e}_3 \cr [{\bf e}_1,{\bf e}_3]={\bf e}_4\cr
[{\bf e}_1,{\bf e}_4]={\bf e}_5 \cr
[{\bf e}_2,{\bf e}_3]={\bf e}_6 \end{matrix}$} &  {\footnotesize
$\left(
\begin{matrix}
1&0&x_2&x_3&x_4& -x_2^2/2\cr
0 &1&0 &0&0&x_3 \cr
0 &0 &1&0&0&0\cr
0 &0 &0 &1&0&0\cr
0 &0 &0 &0 &1&0\cr
0 &0 &0 &0&0 &1
 \end{matrix}
\right)$ }\\

{\footnotesize $A_{6,15}$}&{\footnotesize $\begin{matrix} [{\bf e}_1,{\bf e}_2]={\bf e}_3 + {\bf e}_5\cr
[{\bf e}_1,{\bf e}_3]={\bf e}_4 \hspace*{0.65cm}\cr
[{\bf e}_1,{\bf e}_4]={\bf e}_6 \hspace*{0.65cm}\cr
[{\bf e}_2,{\bf e}_5]={\bf e}_6 \end{matrix}$}&
 {\footnotesize $\left(
\begin{matrix}
1&0&x_2&x_3&x_2&-x_2^2/2 + x_4\cr
0 &1&0 &0&0&x_5 \cr
0 &0 &1&0&0&0\cr
0 &0 &0 &1&0&0\cr
0 &0 &0 &0 &1&0\cr
0 &0 &0 &0&0 &1
 \end{matrix}
\right)$ }\\

{\footnotesize $A_{6,16}$}&{\footnotesize $\begin{matrix} [{\bf e}_1,{\bf e}_3]={\bf e}_4 \cr [{\bf e}_1,{\bf e}_4]={\bf e}_5\cr
[{\bf e}_1,{\bf e}_5]={\bf e}_6 \cr
[{\bf e}_2,{\bf e}_3]={\bf e}_5\cr
[{\bf e}_2,{\bf e}_4]={\bf e}_6 \end{matrix}$} &  {\footnotesize
$\left(
\begin{matrix}
1&0&0&x_3&x_4& x_5\cr
0 &1&0 &0&x_3&x_4 \cr
0 &0 &1&0&0&0\cr
0 &0 &0 &1&0&0\cr
0 &0 &0 &0 &1&0\cr
0 &0 &0 &0&0 &1
 \end{matrix}
\right)$ }\\

{\footnotesize $A_{6,17}$}&{\footnotesize $\begin{matrix}[{\bf e}_1,{\bf e}_2]={\bf e}_3 \cr [{\bf e}_1,{\bf e}_3]={\bf e}_4\cr
[{\bf e}_1,{\bf e}_4]={\bf e}_6 \cr
[{\bf e}_2,{\bf e}_5]={\bf e}_6 \end{matrix}$} & {\footnotesize
$\left(
\begin{matrix}
1&0&x_2&x_3&0& x_4\cr
0 &1&0 &0&0&x_5 \cr
0 &0 &1&0&0&0\cr
0 &0 &0 &1&0&0\cr
0 &0 &0 &0 &1&0\cr
0 &0 &0 &0&0 &1
 \end{matrix}
\right)$ }\\

{\footnotesize $A_{6,18}$}&{\footnotesize $\begin{matrix}[{\bf e}_1,{\bf e}_2]={\bf e}_3 \;\;\;\cr
 [{\bf e}_1,{\bf e}_3]={\bf e}_4\;\;\;\cr
[{\bf e}_1,{\bf e}_4]={\bf e}_6 \;\;\;\cr
[{\bf e}_2,{\bf e}_3]={\bf e}_5 \;\;\;\cr
[{\bf e}_2,{\bf e}_5]=a\;{\bf e}_6 \cr a \neq 0 \end{matrix}$} & {\footnotesize
$\left(
\begin{matrix}
1&0&x_2&x_3&-x_2^2/2& a x_2^3/6 +x_4\cr
0 &1&0 &0&x_3&a x_5 \cr
0 &0 &1&0&0&0\cr
0 &0 &0 &1&0&0\cr
0 &0 &0 &0 &1&0\cr
0 &0 &0 &0&0 &1
 \end{matrix}
\right)$ }\\

\smallskip \\
\hline\hline
\end{tabular}

  \begin{tabular}{l l l}
    \hline\hline
{\footnotesize ${\bf g}$ }&{\footnotesize Non-zero commutation relations } &{\footnotesize vielbeins matrix }\\ \hline

{\footnotesize $A_{6,19}$}&{\footnotesize $\begin{matrix} [{\bf e}_1,{\bf e}_2]={\bf e}_3 \cr [{\bf e}_1,{\bf e}_3]={\bf e}_4\cr
[{\bf e}_1,{\bf e}_4]={\bf e}_5 \cr
[{\bf e}_1,{\bf e}_5]={\bf e}_6\cr
[{\bf e}_2,{\bf e}_3]={\bf e}_6 \end{matrix}$} & {\footnotesize
$\left(
\begin{matrix}
1&0&x_2&x_3&x_4& -x_2^2/2+x_5\cr
0 &1&0 &0&0&x_30 \cr
0 &0 &1&0&0&0\cr
0 &0 &0 &1&0&0\cr
0 &0 &0 &0 &1&0\cr
0 &0 &0 &0&0 &1
 \end{matrix}
\right)$ }\\

{\footnotesize $A_{6,20}$}&{\footnotesize $ \begin{matrix} [{\bf e}_1,{\bf e}_2]={\bf e}_3 \cr [{\bf e}_1,{\bf e}_3]={\bf e}_4\cr
[{\bf e}_1,{\bf e}_4]={\bf e}_5 \cr
[{\bf e}_1,{\bf e}_5]={\bf e}_6 \cr
[{\bf e}_2,{\bf e}_3]={\bf e}_5 \cr
[{\bf e}_2,{\bf e}_4]={\bf e}_6 \end{matrix}$} & {\footnotesize
$\left(
\begin{matrix}
1&0&x_2&x_3&-x_2^2/2+x_4& x_5\cr
0 &1&0 &0&x_3&x_4 \cr
0 &0 &1&0&0&0\cr
0 &0 &0 &1&0&0\cr
0 &0 &0 &0 &1&0\cr
0 &0 &0 &0&0 &1
 \end{matrix}
\right)$ }\\

{\footnotesize $A_{6,24}$}&{\footnotesize $[{\bf e}_1,{\bf e}_2]={\bf e}_3 $} & {\footnotesize
$\left(
\begin{matrix}
1&0&x_2&0&0& 0\cr
0 &1&0 &0&0&0 \cr
0 &0 &1&0&0&0\cr
0 &0 &0 &1&0&0\cr
0 &0 &0 &0 &1&0\cr
0 &0 &0 &0&0 &1
 \end{matrix}
\right)$ }\\

{\footnotesize $A_{6,25}$}&{\footnotesize $ \begin{matrix} [{\bf e}_1,{\bf e}_2]={\bf e}_3 \cr [{\bf e}_4,{\bf e}_5]={\bf e}_6
\end{matrix}$} & {\footnotesize
$\left(
\begin{matrix}
1&0&x_2&0&0& 0\cr
0 &1&0 &0&0&0 \cr
0 &0 &1&0&0&0\cr
0 &0 &0 &1&0&x_5\cr
0 &0 &0 &0 &1&0\cr
0 &0 &0 &0&0 &1
 \end{matrix}
\right)$ }\\

{\footnotesize $A_{6,26}$}&{\footnotesize $\begin{matrix}[{\bf e}_1,{\bf e}_2]={\bf e}_3 \cr [{\bf e}_1,{\bf e}_3]={\bf e}_4
\end{matrix}$} & {\footnotesize
$\left(
\begin{matrix}
1&0&x_2&x_3&0&0\cr
0 &1&0 &0&0&0 \cr
0 &0 &1&0&0&0\cr
0 &0 &0 &1&0&0\cr
0 &0 &0 &0 &1&0\cr
0 &0 &0 &0&0 &1
 \end{matrix}
\right)$ }\\

{\footnotesize $A_{6,27}$}&{\footnotesize $\begin{matrix} [{\bf e}_3,{\bf e}_5]={\bf e}_1 \cr [{\bf e}_4,{\bf e}_5]={\bf e}_2
\end{matrix}$} & {\footnotesize
$\left(
\begin{matrix}
1&0&0&0&0& 0\cr
0 &1&0 &0&0&0 \cr
x_5 &0 &1&0&0&0\cr
0 &x_5 &0 &1&0&0\cr
0 &0 &0 &0 &1&0\cr
0 &0 &0 &0&0 &1
 \end{matrix}
\right)$ }\\

{\footnotesize $A_{6,28}$}&{\footnotesize $\begin{matrix} [{\bf e}_2,{\bf e}_5]={\bf e}_1 \cr [{\bf e}_3,{\bf e}_5]={\bf e}_2\cr
[{\bf e}_4,{\bf e}_5]={\bf e}_3 \end{matrix}$} & {\footnotesize
$\left(
\begin{matrix}
1&0&0&0&0&0\cr
x_5 &1&0 &0&0&0 \cr
x_5^2/2 &x_5 &1&0&0&0\cr
x_5^3/6 &x_5^2/2 &x_5 &1&0&0\cr
0 &0 &0 &0 &1&0\cr
0 &0 &0 &0&0 &1
 \end{matrix}
\right)$ }\\

{\footnotesize $A_{6,31}$}&{\footnotesize $\begin{matrix} [{\bf e}_3,{\bf e}_4]={\bf e}_1 \cr [{\bf e}_2,{\bf e}_5]={\bf e}_1\cr
[{\bf e}_3,{\bf e}_5]={\bf e}_2 \end{matrix}$} & {\footnotesize
$\left(
\begin{matrix}
1&0&0&0&0& 0\cr
x_5 &1&0 &0&0&0 \cr
x_4+x_5^2/2 &x_5 &1&0&0&0\cr
0 &0 &0 &1&0&0\cr
0 &0 &0 &0 &1&0\cr
0 &0 &0 &0&0 &1
 \end{matrix}
\right)$ }\\

{\footnotesize $A_{6,32}$}&{\footnotesize $ \begin{matrix} [{\bf e}_3,{\bf e}_4]={\bf e}_1 \cr [{\bf e}_2,{\bf e}_5]={\bf e}_1\cr
[{\bf e}_3,{\bf e}_5]={\bf e}_2\cr [{\bf e}_4,{\bf e}_5]={\bf e}_3 \end{matrix}$} & {\footnotesize
$\left(
\begin{matrix}
1&0&0&0&0&0\cr
x_5 &1&0 &0&0&0 \cr
x_4+x_5^2/2 &x_5 &1&0&0&0\cr
x_5^3/6 &x_5^2/2 &x_5 &1&0&0\cr
0 &0 &0 &0 &1&0\cr
0 &0 &0 &0&0 &1
 \end{matrix}
\right)$ }\\

\smallskip \\
\hline\hline
\end{tabular}

\end{document}